\documentclass{amsart}
\usepackage{amssymb,latexsym}
\usepackage[mathscr]{eucal}
\usepackage{verbatim}

\newtheorem*{theorem*}{Theorem}
\newtheorem*{corollary*}{Corollary}
\newtheorem{theorem}{Theorem}[section]
\newtheorem{corollary}[theorem]{Corollary}
\newtheorem{lemma}[theorem]{Lemma}
\newtheorem{proposition}[theorem]{Proposition}
\newtheorem{remark}[theorem]{Remark}
\newtheorem{definition}[theorem]{Definition}
\newtheorem{example}[theorem]{Example}

\DeclareMathOperator{\cl}{cl}
\DeclareMathOperator{\diag}{diag}
\DeclareMathOperator{\co*}{c_\mathrm o^*}

\newcommand{\1}{\mspace{2mu}}
\newcommand{\be}{\begin{equation}\label}
\newcommand{\ee}{\end{equation}}
\newcommand{\bq}{\begin{equation*}}
\newcommand{\eq}{\end{equation*}}
\newcommand{\ba}{\begin{align*}}
\newcommand{\ea}{\end{align*}}
\newcommand{\bp}{\begin{proof}}
\newcommand{\ep}{\end{proof}}
\newcommand{\bL}{\begin{lemma}\label}
\newcommand{\eL}{\end{lemma}}
\newcommand{\bP}{\begin{proposition}\label}
\newcommand{\eP}{\end{proposition}}
\newcommand{\bC}{\begin{corollary}\label}
\newcommand{\eC}{\end{corollary}}
\newcommand{\bT}{\begin{theorem}\label}
\newcommand{\eT}{\end{theorem}}
\newcommand{\bR}{\begin{remark}\label}
\newcommand{\eR}{\end{remark}}
\newcommand{\bD}{\begin{definition}\label}
\newcommand{\eD}{\end{definition}}
\newcommand{\bE}{\begin{example}\label}
\newcommand{\eE}{\end{example}}

\begin{document}

\title{Majorization and arithmetic mean ideals}

\author{Victor Kaftal}
\address{University of Cincinnati\\
Department of Mathematics\\
Cincinnati, OH, 45221-0025\\
USA}
\email{victor.kaftal@uc.edu}
\author{Gary Weiss}
\email{gary.weiss@uc.edu}

\keywords{Majorization of sequences, operator ideals, arithmetic mean}
\subjclass{Primary: 15A51, 47L20 }
\date{\today}
\begin{abstract}
Following \textit{An infinite dimensional Schur-Horn theorem and majorization theory} \cite{vKgW10-Maj},
this paper further studies majorization for infinite sequences.
It extends to the infinite case classical results on ``intermediate sequences" for finite sequence majorization.
These and other infinite majorization properties are then linked to notions of infinite convexity and invariance properties under various classes of substochastic matrices to characterize arithmetic mean closed operator ideals and arithmetic mean at infinity closed operator ideals.
\end{abstract}

\maketitle


\section{Introduction}\label{S: 1}
Operator ideals - two sided ideals of the algebra $B(H)$ of bounded linear operators on an infinite dimensional separable Hilbert space $H$ - have been the object of much study in the last 70 years. A central role in the development of their theory has always been played by the study of their commutator spaces
(see listed chronologically \cite{pH54}, \cite {BP65}, \cite{PT71}, \cite{jA77}-\cite{AV86}, \cite{gW80}-\cite{gW86}, \cite{nK89}).
The introduction in the 1980's of cyclic cohomology by
A. Connes (see \cite{aC82}-\cite {aC85}) and the connection in the 1990's with algebraic K-theory by M. Wodzicki (see \cite {mW94}) provided further motivation for this study. The full characterization of the commutator space of an ideal was achieved by Dykema, Figiel, Weiss, and Wodzicki in \cite {DFWW} in terms of arithmetic means. In particular, it lead to the notion of \emph{arithmetic mean closed} ideals
(see Section \ref{S: 2} preceding (\ref{e: 4}))
and numerous operator ideals studied by many in the past were shown to possess this property.

Although operator ideals arise in a quintessentially \emph{non-commutative} setting, it evolved from the work of Calkin \cite{jC41} and Schatten \cite{rS60} that each ideal is characterized by a certain cone of sequences that are monotone, positive, real-valued and tend to 0,
the \emph{characteristic set} of the ideal. (See \cite[Section 2]{DFWW} for background and development.)
It turns out that the property for an ideal of being arithmetic mean closed translates seamlessly into the hereditariness (i.e., solidity) of its characteristic set under \emph {majorization} (see Definition \ref{D:1.2} below).

Recall that if $\xi$ and $ \eta$ are two finite real valued sequences, i.e., two vectors in $ \mathbb R^N$ then $\xi$ is said to be majorized by $\eta$,
(in symbols $\xi\preccurlyeq \eta$), if $ \sum_{j=1}^n \xi^*_j \le \sum_{j=1}^n \eta^*_j $ for $1\le n\le N$ and $ \sum_{j=1}^N \xi_j = \sum_{j=1}^N \eta_j $,
where $\xi^*, \eta^*$ denote the monotone nonincreasing rearrangement of $\xi, \eta$.
The theory of majorization started more than a century ago motivated by wealth distribution (Lorenz \cite {Lm1905}), inequalities involving convex functions and the study of stochastic matrices (Hardy, Littlewood and P\'olya \cite {HLP52}), convex combinations of permutation matrices (Birkhoff \cite{Bg46}), and the relation between diagonals of selfadjoint matrices and their eigenvalue lists (Schur \cite {Si23} and Horn \cite{Horn},
and others \cite {Mi58}-\cite{Mi60}, \cite {aM64}, \cite {GiMa64}).
In more recent years majorization theory and also its extension to infinite sequences has turned up prominently in several areas of operator theory and frame theory
(\cite {Kr02a}-\cite {Kr02b}, \cite {AK02}, \cite {CL02}, \cite {AMRS}, \cite{KkLd04}).
We contributed to this study in \cite {vKgW10-Maj}, where we extended many classical facts on (finite) majorization to infinite sequences decreasing to zero and used these extensions to prove an infinite dimensional Schur-Horn theorem.

One of our initial interests in developing infinite majorization theory was to obtain tools for the study of arithmetic mean ideals, one goal of the present paper which is part of a long-term study of arithmetic mean ideals and arithmetic mean at infinity ideals \cite{vKgW02}-\cite{vKgW07-OT21}.

This paper is organized as follows.
In Section \ref {S: 2} we present background and notations, including basic information on operator ideals and on infinite majorization, ``majorization at infinity," strong versions of both, and the yet stronger notion of block majorization.

In Section \ref {S: 3} we extend to infinite majorization several classical results in finite majorization theory on ``intermediate sequences."
Intermediate sequences deepen our understanding of majorization theory by relating various of its properties.
Techniques based on finite intermediate sequences have been applied in an operator theory context also in the work of Arveson and Kadison \cite{AK02}.

In Section \ref {S: 4} we employ techniques based on infinite intermediate sequences to characterize arithmetic mean closed ideals and arithmetic mean at infinity closed ideals in terms of the diagonals of their positive elements (Propositions \ref {P:4.1}, \ref{P:4.3}).

As a consequence, we obtain that an ideal is arithmetic mean closed if and only if the ideal is diagonally invariant, that is,
for every element of the ideal, the main diagonal of that element (for any fixed basis) is itself in the ideal (Theorem \ref {T:4.5}).

A further characterization is that an ideal is arithmetic mean closed if and only if it is invariant under the action of all the substochastic matrices (Theorem \ref {T:4.7}).

We analyze also invariance under smaller classes of substochastic matrices.
For instance, invariance under the action of orthostochastic matrices (Definition \ref {D:stoch}) is also necessary and sufficient.
Invariance under the action of infinite convex combinations of permutations (a class studied by Birkhoff in the finite case, see remarks after Theorem \ref{T:4.7}) is in general not sufficient,
but is sufficient for soft ideals (Definition \ref {D:soft}, Theorem \ref {T:4.9}).
Soft ideals form a class of ideals that includes many of the classical operator ideals and has been studied in our work in \cite {vKgW04-Soft}.

\section{Background and Notations}\label{S: 2}

We start by setting notations for majorization of infinite sequences.
These are natural extensions of the notions mentioned in the Introduction for finite sequences and are the same notations that we used in \cite {vKgW10-Maj}.
We found that the terminology and notations describing majorization for finite sequences vary considerably in the literature.

To avoid always having to pass to monotone rearrangements, herein we will focus on sequences decreasing monotonically to 0 and will denote by $\co*$ their positive cone and by $(\ell^1)^*$ the subcone of summable decreasing sequences. (We note explicitly that $\co*$ and $(\ell^1)^*$ do not denote herein the duals of $\text{c}_{\text{o}}$ and $\ell^1$.)

\bD{D:1.2}
For $\xi, \eta \in \co*$ we say that
\begin{itemize}
\item
$\eta$ \textit{majorizes} $\xi$ ($\xi \prec \eta$) if $ \, \sum_{j=1}^n \xi_j \le \sum_{j=1}^n \eta_j $ for every $n\in \mathbb N; $
\item
$\eta$ \textit{strongly majorizes} $\xi$ ($\xi\preccurlyeq \eta)
\text{ if } \xi \prec \eta \, \text { and } {\varliminf }\sum_{j=1}^n ( \eta_j -\xi_j) = 0$;
\item $\eta$ \textit{block majorizes} $\xi$ ($\xi\prec_b \eta)
\text { if } \xi \prec \eta \, \text { and } \sum_{j=1}^{n_k} \xi_j = \sum_{j=1}^{n_k} \eta_j$
for some sequence $ \mathbb N \ni n_k \uparrow \infty$.
\end{itemize}
\noindent For $\xi, \eta \in (\ell^1)^*$ we say that
\begin{itemize}
\item $\eta$ \textit{majorizes at infinity} $\xi$ ($ \xi \prec_\infty\eta$)
if $ \sum_{j=n}^\infty \xi_j \le \sum_{j=n}^\infty \eta_j$ for every $ n\in \mathbb N;$
\item $\eta$ \textit{strongly majorizes at infinity} $\xi$ ($ \xi \preccurlyeq_\infty \eta$)
if $\xi \prec_\infty\eta $ and
$ \sum_{j=1}^\infty \xi_j = \sum_{j=1}^\infty\eta_j$.
\end{itemize}
If $\xi, \eta \in \text{c}_{\text{o}}^+ $ (resp. $ (\ell^1)^+$), then we say that any of the above relations hold for $\xi$ and $\eta$ whenever they hold for their monotone rearrangements $\xi^*$ and $\eta^*$.
\eD
Immediate consequences of Definition \ref {D:1.2} are:
\begin{align}
&\text{If $\xi, \eta \in \co*$, then}\quad \xi\prec_b \eta \Rightarrow \xi\preccurlyeq \eta \Rightarrow \xi \prec \eta.\label{e: 1}\\
&\text{If $\xi, \eta \in (\ell^1)^*$, then } \xi\prec \eta \text{ and }\sum_{j=1}^\infty \xi_j = \sum_{j=1}^\infty\eta_j\Leftrightarrow \xi\preccurlyeq \eta\, \Leftrightarrow \, \eta \preccurlyeq_\infty\xi.\label{e: 2}
\end{align}

For nonsummable monotone decreasing sequences, the condition
${\varliminf}\sum_{j=1}^n (\eta_j - \xi_j) = 0$
retains many of the key properties of ``equality at the end" for finite and for summable sequences and plays a major role in \cite{vKgW10-Maj}.
There we showed that if $\xi, \eta\in \co*$ and $\xi\prec \eta$, then $\xi = Q(\xi,\eta)\eta$ for some infinite matrix $Q(\xi,\eta)$,
the Schur-square of a real-entry canonical co-isometry $W(\xi, \eta)$ (i.e., $Q(\xi, \eta)_{ij}= \big(W(\xi, \eta)_{ij}\big)^2$).
This representation is directly related to diagonals of operators in unitary orbits (see Lemma \ref{L:2.3} and the succeeding comment).

Majorization at infinity, aka ``tail majorization," was first introduced and studied for finite sequences and
appears in \cite{AGPS87}, \cite {nK89}, \cite{DFWW}, \cite {mW02} among others.
Here it provides the natural characterization of the notion of arithmetic mean at infinity closure for operator ideals contained in the trace class (see Proposition \ref {P:4.3} and \cite{vKgW02}-\cite {vKgW07-OT21}).

Block-majorization is a natural way to bring to bear the results of finite majorization theory on infinite sequences
and it also arises naturally in the study in \cite[Section 4]{vKgW10-Maj} of the canonical co-isometry $W(\xi, \eta)$ associated with two sequences $\xi\prec \eta$.
It also plays a key role in the present paper so is further developed in Section \ref{S: 3}.

We need to review here some of the connections between majorization, stochastic matrices and conditional expectations.
Fix an othonormal basis $\{e_n\}$ for the Hilbert space $H$.
\bD{D:stoch}
A matrix $P$ with nonnegative entries is called
\begin{itemize}
\item substochastic if all its row and column sums are bounded by 1;
\item column-stochastic if it is substochastic with all column sums equal to 1;
\item row-stochastic if it is substochastic with all row sums equal to 1;
\item doubly stochastic if it is both column and row-stochastic;
\item block stochastic if it is the direct sum of finite doubly stochastic matrices.
\item orthostochastic if $P_{ij}= (U_{ij})^2$ for some orthogonal matrix $U$ (i.e., a unitary matrix with real entries).
\end{itemize}
\eD
\noindent So any orthostochastic matrix is also doubly stochastic.\\

Notice that we can apply a substochastic matrix $P$ to any sequence $\rho \in \ell^\infty$,
where by $P\rho$ we simply mean the sequence $\big<\sum_{j=1}^\infty P_{ij}\rho_j\big>_i$.
If $\rho\in \text{c}_{\text{o}}$ and $P$ is substochastic, then $P\rho \in \text{c}_{\text{o}}$.
By Schur's test (e.g., see \cite [Problem 45]{pH82}) substochastic matrices viewed as Hilbert space operators are contractions.

\bL {L:Mark}\cite[Lemma 3.1]{aM64}.
For $\xi, \eta\in \co*$,
$$\xi \prec \eta \text{\quad if and only if \quad there is a substochastic matrix} ~P~ \text{for which}\quad \xi = P\eta.$$
\eL

In fact, more
can be said about this regarding the important class of orthostochastic
matrices which are at the core of the celebrated
Schur-Horn Theorem and its infinite dimensional extension in \cite{vKgW10-Maj}.

\bT{T:SH}\cite [Theorems 3.9 \& 5.3]{vKgW10-Maj}.
Let $\xi, \eta\in \co*$ and assume that $\xi \prec \eta$.
\item[(i)] If $\xi\not\in \ell^1$, there is an orthostochastic matrix $Q$ such that $\xi = Q\eta.$
\item[(ii)] If $\xi\in \ell^1$, there is an orthostochastic matrix $Q$ such that $\xi = Q\eta$ if and only if $\xi \preccurlyeq \eta.$
\eT

Notice also an immediate consequence of the (finite) Schur-Horn Theorem \cite [Theorem 4] {Horn}:
\begin{align}\label{e: 3}
\text{For }\xi, \eta \in \co*,\quad \xi\prec_b \eta & \Leftrightarrow \xi = Q\eta \text{ for a block orthostochastic matrix }Q\\
& \Leftrightarrow \xi = Q\eta \text{ for a block stochastic matrix }Q.\notag
\end{align}

\bL{L:2.10}\cite[Lemma 2.10]{vKgW10-Maj}. Let $\xi, \eta \in \co*$ with $\xi = P\eta$ for some column-stochastic matrix $P$. If $\xi \in (\ell^1)^*$ (resp. $\eta\in(\ell^1)^*$), then $\eta \in (\ell^1)^*$ (resp. $\xi \in (\ell^1)^*$) and $\xi \preccurlyeq \eta$.
\eL

Furthermore, denote by $E$ the conditional expectation on the maximal abelian algebra of diagonal operators (relative to the fixed orthonormal basis $\{e_n\}$).
$E$ is the operation of ``taking the main diagonal" in the matrix representation relative to $\{e_n\}$.
It is elementary to verify that

\bL{L:2.3} \cite[Lemmas 2.4, 2.3]{vKgW10-Maj}.
Let $L\in B(H)$ be a contraction and let $Q$ be its Schur-square, i.e., $Q_{ij} = |L_{ij}|^2$ for all $i,j$.
\item [(i)] $Q$ is substochastic.
\item [(ii)] Respectively, $Q$ is column-stochastic, row-stochastic, block orthostochastic, orthostochastic if and only if, respectively, $L$ is an isometry, co-isometry, a direct sum of finite orthogonal matrices, an orthogonal matrix.
\item [(iii)] For $\xi, \eta \in \co*$, $\diag \xi = E(L\diag\eta \, L^*)$ if and only if $\xi = Q \eta$.
\eL
\medskip
Next we recall some operator ideal facts.
Calkin \cite{jC41} established a correspondence between the two-sided
ideals of $B(H)$ for a separable infinite-dimensional
Hilbert space $H$ and the \textit{characteristic sets}.
These are the positive cones of $\text{c}_{\text{o}}^*$
that are hereditary (i.e., solid) and invariant under \textit{ampliations}, i.e., the maps
\[
\text{c}_{\text{o}}^* \owns \xi \rightarrow
D_m\xi:=~<\xi_1,\ldots,\xi_1,\xi_2,\ldots,\xi_2,\xi_3,\ldots,\xi_3,\ldots>
\]
where each entry $\xi_i$ of $\xi$ is repeated $m$-times.
The order-preserving lattice isomorphism $I \rightarrow \Sigma(I)$ maps each ideal to its characteristic set $\Sigma(I) := \{s(X) \mid X \in I\}$,
where $s(X)$ denotes the singular number sequence of $X$, namely the sequence of the eigenvalues of $|X|=(X^*X)^{1/2}$, in decreasing order, and repeated according to multiplicity.
For a fixed orthonormal basis, the inverse of this lattice isomorphism is the map $\Sigma\rightarrow I(\Sigma)$ where
$I(\Sigma)$ is the ideal generated by $\{\diag \xi \mid \xi \in \Sigma\}$, where $\diag \xi$ is the diagonal operator
(with respect to that basis) with diagonal $\xi$.

Two sequence operations, the arithmetic mean restricted to $\co*$ and the
arithmetic mean at infinity restricted to $(\ell^1)^*$, respectively
\[
\xi_a := \langle\frac{1}{n}\sum_{j=1}^n \xi_j\rangle
\quad \text{and} \quad \xi_{a_\infty} :=
\langle\frac{1}{n}\sum_{j=n+1}^\infty \xi_j\rangle,
\]
are essential for the characterization of commutator ideals (also called commutator spaces) (see \cite{DFWW})
and for defining next the consequent natural classes of arithmetic mean ideals and arithmetic mean at infinity ideals
(see also \cite{vKgW02}-\cite{vKgW07-OT21}).

If $I$ is an ideal, then the arithmetic mean ideals $_aI$ and $I_a$,
the \textit{pre-arithmetic mean ideal} and \textit{arithmetic mean ideal}
of $I$, are the ideals with characteristic sets
\[
\Sigma(_aI) := \{\xi \in \text{c}_{\text{o}}^* \mid \xi_a \in \Sigma(I)\},
\]
\[
\Sigma(I_a) := \{\xi \in \text{c}_{\text{o}}^* \mid \xi =
O(\eta_a)~\text{for some}~ \eta \in \Sigma(I)\},
\]
where the notation $\xi=O(\zeta)$ means $\xi_n\le M\zeta_n$ for some $M>0$ and all $n$.
The ideal $I^-:=\,_a(I_a)$ is called the \textit{am-closure} of $I$ and an ideal $I$ is called
\textit{am-closed} (short for arithmetic mean closed) if $I = I^-$.
(There is are analogous notions of \textit{am-interior} and \textit{am-open} ideals [ibid].)
It is easy to see that $I\subset I^-$ and that the am-closure of the finite rank ideal $F$ is the trace-class ideal $\mathscr L_1$ (see \cite{vKgW04-Traces}).
Thus
\be{e: 4}
\mathscr L_1 = F^- \quad\text{is the smallest am-closed ideal.}
\ee
Many of the ideals arising from classical sequence spaces are am-closed and am-closed ideals play a central role in the study of single commutators in operator ideals, e.g., \cite[Theorem 7.3 and Corollary 7.10]{DFWW}.

Similarly, the arithmetic mean at infinity ideals $_{a_\infty}I$ and $I_{a_\infty}$
are the ideals with characteristic sets
\[
\Sigma(_{a_\infty}I) := \{\xi \in (\ell^1)^* \mid \xi_{a_\infty} \in
\Sigma(I)\}
\]
\[
\Sigma(I_{a_\infty}) := \{\xi \in \text{c}_{\text{o}}^{*} \mid \xi =
O(\eta_{a_\infty}) ~\text{for some}~ \eta \in \Sigma(I \cap \mathscr L_1)\}.
\]
The am-$\infty$ closure of $I$ is
$I^{-\infty} :=\, _{a_\infty}(I_{a_\infty})$ and an ideal $I$ is called am-$\infty$ closed if $I = I^{-\infty}$.
(There is an analogous notion of am-$\infty$ interior and am-$\infty$ open ideals.)
The map $I\rightarrow I^{-\infty}$ is also idempotent, but by definition,
$I^{-\infty}\subset \mathscr L_1$ and hence we only have $I\cap \mathscr L_1 \subset I^{-\infty}$.
For the definition and basic properties of arithmetic mean at infinity ideals see \cite{vKgW04-Soft} and \cite {vKgW04-Traces}.
\cite {vKgW02}-\cite {vKgW07-OT21} develop the theory of arithmetic mean and arithmetic mean at infinity ideals.

The connection between the lattice of arithmetic mean ideals (resp., arithmetic mean at infinity ideals) and majorization (resp., majorization at infinity) follows from the obvious relations:
\begin{align}
\phantom{aaaaaa}&\text{If }\xi, \eta \in \co* &&\text{then } \xi \prec \eta &&\text {if and only if} &\xi_a&\le \eta_a;\phantom{aaaaaaaaaaaaaaaaaaaaaaaaaaaa}\label{e: 5}\\
\phantom{aaaaaa}&\text{if }\xi, \eta \in (\ell^1)^* &&\text{then } \xi \prec_\infty \eta &&\text {if and only if} &\xi_{a_\infty}&\le \eta_{a_\infty}.\label{e: 6}
\end{align}
The consequences of this connection will be analyzed in Section \ref {S: 4}.

Finally, let us review another class of operator ideals, the \textit{soft} ideals, which are often encountered in the study of arithmetic mean ideals.
In this paper they will play an important role in the connection between arithmetic mean ideals and diagonal invariance of ideals (see Theorem \ref{T:4.5}) under some classes of substochastic matrices (see Definition \ref{D:4.6}, Lemma \ref{P:4.8} and Theorem \ref{T:4.9}).

Soft ideals, that is \textit{soft-edged} and \textit{soft-complemented} ideals, are defined as follows:

\bD{D:soft}\cite{vKgW04-Soft}.\\
An ideal $I$ is \textit{soft-edged} if every $\xi\in \Sigma(I)$ is a product $\xi = \alpha \eta$ for some $\alpha \in \co*$ and some $\eta \in \Sigma(I)$.\\
An ideal $I$ is \textit{soft-complemented} if $\xi\in \Sigma(I)$ if and only if $ \alpha\xi \in \Sigma(I)$ for all $\alpha \in \co*$.
\eD
\noindent The soft-complemented terminology derives from its natural equivalence:
If $\xi \in \co*$, then $\xi\not\in \Sigma (I)$ if and only if $\alpha\xi \not \in\Sigma (I)$ for some $\alpha \in \co*$.
Soft ideals have historically played an important role in the theory of operator ideals, albeit that role was mainly implicit.
In \cite{vKgW04-Soft} we proved that many classical sequence ideals are soft,
including countably generated ideals, minimal Banach ideals
$\mathfrak{S}_\phi^{(o)}$ and maximal Banach ideals $\mathfrak S_{\phi}$
for symmetric norming functions $\phi$, Orlicz ideals $\mathscr L_M$
and small Orlicz ideals $\mathscr L_M^{(o)}$ for monotone
nondecreasing functions $M$ on $[0,\infty)$ with $M(0) = 0$,
Lorentz ideals $\mathscr L(\phi)$ for monotone nondecreasing nonnegative sequences $\phi$ satisfying the $\Delta_2$-condition,
Marcinkiewicz ideals, K\"{o}the duals and idempotent ideals.
For further properties of soft-edged and soft-complemented ideals we refer the reader to \cite {vKgW04-Soft}. (See also \cite {vKgW02}, \cite{vKgW04-Traces}.)

\section{Intermediate sequences} \label{S: 3}

The main goal of this section is to extend to infinite sequences the finite intermediate sequence results we list in the following proposition.
(See also \cite[5.A.9 and 5.A.9a]{MO79}.)
\bP{P:3.1} Let $\xi, \eta \in \mathbb (R^N)^+$ be two nonincreasing finite sequences.
\item[(A)] If $\xi \prec \eta$, then there are finite nonincreasing sequences $\zeta, \rho \in \mathbb (R^N)^+$ with
\item[(i)] $\xi \le \rho\preccurlyeq \eta$
\item[(ii)] $\xi\preccurlyeq \zeta \le \eta$
\item[(B)] If $\xi \prec_\infty\eta$, then there are finite nonincreasing sequences $\zeta, \rho \in \mathbb (R^N)^+$ with
\item[(i)] $\xi \le \rho \preccurlyeq_\infty \eta$ \quad(equivalently, $\xi \le \rho$ and $\eta\preccurlyeq \rho$)
\item[(ii)] $\xi \preccurlyeq_\infty \zeta \le \eta$ \quad(equivalently, $\zeta\preccurlyeq \xi$ and $\zeta \le \eta$)
\eP
\noindent (A)(i) is due to Fan \cite {Fa51} (see also \cite [5.A.9]{MO79}).\\
(A)(ii) is due to Mirsky \cite {Mi60} and the one-line proof given for $\mathbb R^N$ in \cite [5.A.9]{MO79} easily extends to $\mathbb (R^N)^+$ (see also \cite [Lemma 4.3]{AK02}) and to $(\ell^1)^*$. \\
(B)(i)-(ii) are stated without proof or attribution in \cite [5.A.9a]{MO79}.
The proof of (B)(i) is immediate by choosing $\rho:=\,< \xi_1 + \sum_{j=1}^N(\eta_j-\xi_j), \xi_2, \ldots,\xi_N>$.
The proof of (B)(ii) is similar to that of (A)(i) but since we use (B)(ii) in Lemma \ref {L:3.10.0} and Theorem \ref {T:3.11}, for completeness's sake and the readers' convenience, we prove it here.
\bp[Proof of (B)(ii)] If $\sum_{j=1}^N(\eta_j-\xi_j) = 0$ choose $\zeta= \eta$.
So we can assume without loss of generality that $\sum_{j=1}^N(\eta_j-\xi_j) > 0$.
The proof is by induction on $N$. If $N=1$ then $\xi_1 < \eta_1$, thus $\zeta:=\xi$ satisfies the condition.
Assume the property true for all $k \le N-1$ for some $N > 1$ and let $ \alpha:= \underset{1\le n \le N }{\min} \sum_{j=n}^N(\eta_j-\xi_j)$.
In particular, $\eta_N -\alpha\ge \xi_N \ge 0$.
Passing if necessary to $\tilde{\eta}:=~ <\eta_1, \eta_2, \ldots, \eta_{N-1}, \eta_N-\alpha>$, and since
$\xi\prec_\infty \tilde{\eta}\le \eta$ and $\underset{1\le n \le N}{\min} \sum_{j=n}^N(\tilde{\eta}_j-\xi_j)=0$, we can assume without loss of generality that $\alpha = 0$.
Let $m$ be an integer for which the minimum $\alpha=0$ is attained, i.e., $ \sum_{j=m}^N(\eta_j-\xi_j)=0$.
Then $1<m\le N$ since by assumption $\sum_{j=1}^N(\eta_j-\xi_j) > 0$.
For every $1\le n \le m-1$ we have
$$\sum_{j=n}^{m-1}(\eta_j-\xi_j) = \sum_{j=n}^N(\eta_j-\xi_j)\ge 0,$$ i.e., $ <\xi_1, \xi_2, \ldots, \xi_{m-1}>~ \prec_\infty ~<\eta_1, \eta_2, \ldots, \eta_{m-1}>$.
By the induction hypothesis there is a monotone nonincreasing sequence $<\zeta_1, \zeta_2, \ldots, \zeta_{m-1}>$ for which
\[
<\xi_1, \xi_2, \ldots, \xi_{m-1}> \;\preccurlyeq_\infty \; <\zeta_1, \zeta_2, \ldots, \zeta_{m-1}> \;\le \; <\eta_1, \eta_2, \ldots, \eta_{m-1}>.
\]
In particular, $\sum_{j=n}^{m-1}(\zeta_j-\xi_j) \ge 0$ for $1\le n \le m-1$, with equality for $n=1$. Now define
\[
\zeta:=\,<\zeta_1, \zeta_2, \ldots, \zeta_{m-1}, \eta_m, \eta_{m+1}, \ldots, \eta_N>.
\]
Then $\zeta \le \eta$,
$$
\sum_{j=n}^N(\zeta_j-\xi_j) =
\begin{cases}
\sum_{j=n}^{m-1}(\zeta_j-\xi_j) \ge 0 \quad &\text{for $1\le n \le m-1$}\\
\sum_{j=n}^N (\eta_j-\xi_j) \ge 0 & \text{for $m\le n \le N$}
\end{cases},
$$
and thus $\xi \prec_\infty \zeta$. Furthermore, $ \sum_{j=1}^N(\zeta_j-\xi_j)= \sum_{j=1}^{m-1}(\zeta_j-\xi_j) = 0$ and so $\xi \preccurlyeq_\infty \zeta$.
It remains to prove that $\zeta$ is nonincreasing, and this requires merely proving that $\zeta_{m-1} \ge \zeta_m$. Since $ \sum_{j=n}^{m-1}(\zeta_j-\xi_j)\ge 0$ for all $n \le m-1$, it follows that $\zeta_{m-1}\ge \xi_{m-1} \ge \xi_m$.
Since $\sum_{j=m}^{N}(\zeta_j-\xi_j)= 0$ while $\sum_{j=m+1}^{N}(\zeta_j-\xi_j)\ge 0$, it follows that $\xi_m \ge \zeta_m$, whence $\zeta_{m-1}\ge \zeta_m$ which completes the proof.
\ep

The first step in the extension of (A)(i) to infinite sequences is the following proposition.

For sequence $\xi$ and integer $n \ge 0$, denote its truncated sequence by $\xi^{(n)} := ~<\xi_{n+1}, \xi_{n+2}, \ldots>$, that is, $\xi^{(n)}_j = \xi_{n+j}$ for $j \ge 1$.

\bP{P:3.2}
If $\xi, \eta \in \co*$ and $\xi \prec \eta$, then the following conditions are equivalent.
\item[(i)] There is a strictly increasing sequence of integers $ n_k$ starting with $n_o=0$ for which
$\xi^{(n_k)} \prec \eta^{(n_k)}.$
\item[(ii)] $\{ \sum_{j=1}^n(\eta_j-\xi_j)\} _{ n \ge m}$ attains a minimum for every $m\in \mathbb N$.
\item[(iii)] $\xi \le \rho\prec_b \eta$ for some $\rho \in \co*$.
\eP
\bp
\item[(i)] $\Leftrightarrow $ (ii)
Notice first that if for some $m\in \mathbb N$, $\xi^{(m)}\prec \eta^{(m)}$, then by definition $\sum_{j=m+1}^n(\eta_j -\xi_j)\ge 0$ for every $n \ge m+1$, hence
$\sum_{j=1}^n (\eta_j -\xi_j) \ge \sum_{j=1}^m (\eta_j -\xi_j) $ for every $n \ge m$. Thus $\{ \sum_{j=1}^n(\eta_j-\xi_j)\} _{ n \ge m}$ attains its minimum for $n=m$.
Thus if (i) holds, then for every $m \in \mathbb N$, $\{ \sum_{j=1}^n(\eta_j-\xi_j)\} _{ n \ge m}$ attains its minimum because for every $k$,
if $n_{k-1}< m\le n_k$, then
\[
\inf_{n \ge m} \sum_{j=1}^n(\eta_j-\xi_j) = \min~\{\min_{m\le n< n_k} \sum_{j=1}^n(\eta_j-\xi_j), \inf_{n\ge n_k}\sum_{j=1}^n(\eta_j-\xi_j)\} \\
= \min_{m\le n\le n_k} \sum_{j=1}^n(\eta_j-\xi_j).
\]
Conversely, if (ii) holds, then set $n_o=0$ and define $n_k$ recursively by setting $n_{k+1}$ to be an index for which
the minimum of $\{ \sum_{j=1}^n(\eta_j-\xi_j)\}_{ n > n_k}$ is attained.
Then $\sum_{j=1}^n (\eta_j -\xi_j) \ge \sum_{j=1}^{n_{k+1}} (\eta_j -\xi_j) $ for all $n > n_k$ and in particular for all $n > n_{k+1}$.
Subtracting obtains
$\sum_{j=n_{k+1}+1}^n (\eta_j -\xi_j) \ge 0$, i.e, $\xi^{(n_{k+1})}\prec \eta^{(n_{k+1})}$.
\item[(iii)] $\Rightarrow $ (i) Let $\sum_{j=1}^{n_k} \rho_j = \sum_{j=1}^{n_k} \eta_j$ for a sequence of indices $\mathbb N \ni n_k \uparrow \infty$ and set $n_0=0$.
For $k=0$, $\xi^{(0)} = \xi \prec \eta = \eta^{(0)}$ and for all $m,k \in \mathbb N$,
\begin{alignat*}{2}
\sum_{j=1}^m\rho^{(n_k)}_j&= \sum_{j=n_k+1}^{n_k+m}\rho_j &&(\text{by definition})\\
&= \sum_{j=1}^{n_k+m}\rho_j -\sum_{j=1}^{n_k}\rho_j \\
&\le \sum_{j=1}^{n_k+m}\eta_j -\sum_{j=1}^{n_k}\rho_j &&(\text{since } \rho\prec \eta)\\
&= \sum_{j=1}^{n_k+m}\eta_j -\sum_{j=1}^{n_k}\eta_j\qquad \qquad &&(\text{by hypothesis})\\
&= \sum_{j=n_k+1}^{n_k+m}\eta_j\\
&= \sum_{j=1}^m\eta^{(n_k)}_j &&(\text{by definition})
\end{alignat*}
Thus $\rho^{(n_k)} \prec \eta^{(n_k)}$.
Since $\xi\le \rho$ one has $\xi^{(n_k)}\le \rho^{(n_k)}$ and hence $\xi^{(n_k)}\prec \eta^{(n_k)}$ by transitivity.

\item (i) $\Rightarrow $ (iii) Applying Proposition \ref{P:3.1}(A)(i) to each pair of finite sequences
\[
<\xi_{n_k+1}, \xi_{n_k+2}, \ldots, \xi_{n_{k+1}}>\, ~\prec~ \,<\eta_{n_k+1}, \eta_{n_k+2}, \ldots, \eta_{n_{k+1}}>
\]
we find a nonincreasing finite sequence $<\rho_{n_k+1}, \rho_{n_k+2}, \ldots, \rho_{n_{k+1}}>
$ for which
\[
<\xi_{n_k+1}, \xi_{n_k+2}, \ldots, \xi _{n_{k+1}}> ~ \le ~ <\rho_{n_k+1}, \rho_{n_k+2}, \ldots, \rho_{n_{k+1}}> ~\preccurlyeq~ <\eta_{n_k+1}, \eta_{n_k+2}, \ldots, \eta_{n_{k+1}}>.
\]
Let $\rho$ be the sequence obtained by splicing together these finite sequences.
By construction, $ \xi \le \rho $ and it is immediate to see that $\rho \prec_b \eta$. Since
$
\rho_{n_k} - \eta_{n_k} = \sum_{j=n_{k-1}+1}^{n_k-1}(\eta_j-\rho_j)\ge 0
$
and hence $$
\rho_{n_k} \ge \eta_{n_k} \ge \eta_{n_k+1} \ge \rho_{n_k+1},
$$
we see that $\rho$ is monotone nonincreasing, thus concluding the proof.
\ep

\bR{R:3.3}
A sufficient condition for (ii) in Proposition \ref {P:3.2} to hold is that $\xi\prec \eta$ and $\sum_{j=1}^{\infty}(\eta_j-\xi_j)=\infty$. \eR

\bT{T: 3.4} If $\xi, \eta \in \co*$ and $\xi \prec \eta$, then there is some $\rho\in \co*$ such that $\xi\le\rho \preccurlyeq \eta$.
\eT

\bp
Let $\alpha:= \varliminf_n\sum_{j=1}^n (\eta_j-\xi_j)$.
If $\alpha=0$ then $\xi\preccurlyeq \eta$, in which case it suffices to choose $\rho=\xi$. \linebreak
If condition (ii) of Proposition \ref {P:3.2} is satisfied, which by Remark \ref {R:3.3} includes the case $\alpha=\infty$,
then there is a $\rho \in \co*$ for which $\xi \le \rho\prec_b \eta$ and hence $\xi\le\rho \preccurlyeq \eta$.
Thus it suffices to assume that $0< \alpha < \infty$ and that condition (ii) of Proposition \ref {P:3.2} fails.
So then there is an $m\in \mathbb N$ for which $\{ \sum_1^n(\eta_j-\xi_j)\} _{ n \ge m}$ has no minimum from which
it follows that $\alpha = \inf _{n\ge m}\sum_{j=1}^n ( \eta_j -\xi_j)$ and $\sum_{j=1}^n ( \eta_j -\xi_j)>\alpha$ for every $n\ge m$.
Let $M$ be the smallest such integer.

If $M=1$, namely $\sum_{j=1}^n ( \eta_j -\xi_j) >\alpha$ for every $n\ge 1$ and $\inf_n\sum_{j=1}^n ( \eta_j -\xi_j) =\alpha$, then let
$$ \rho:=\,< \xi_1+ \alpha , \xi_2, \xi_3, \ldots>.$$
It is clear that $\rho$ is monotone nonincreasing and it is immediate to verify that it satisfies the requirement $\xi\le\rho \preccurlyeq \eta$.

If $M > 1$, then $\sum_{j=1}^{M-1}(\eta_j-\xi_j) \le \alpha$ since otherwise $\sum_{j=1}^{M-1}(\eta_j-\xi_j) > \alpha$, whence $\{ \sum_{j=1}^n(\eta_j-\xi_j)\} _{ n \ge M-1}$ also has no minimum, which contradicts the minimality of $M$. Thus
\be{e: 7}
\sum_{j=1}^{M-1}(\eta_j-\xi_j) \le \alpha < \sum_{j=1}^n(\eta_j-\xi_j)\quad \text{ for all } n \ge M.
\ee
Then by Proposition \ref {P:3.1}(A)(i) applied to the sequences
$
<\xi_1, \xi_2, \ldots, \xi_{M-1}>\, \prec \,<\eta_1, \eta_2, \ldots, \eta_{M-1}>,
$
we can find $<\rho_1, \rho_2, \ldots, \rho_{M-1}> ~\in \mathbb (R^{M-1})^+$ with
\be{e: 8}
<\xi_1, \xi_2, \ldots, \xi_{M-1}>\,\le \,<\rho_1, \rho_2, \ldots, \rho_{M-1}> \,\preccurlyeq\,<\eta_1, \eta_2, \ldots, \eta_{M-1}>.
\ee
Set
$$\rho:=\,<\rho_1, \rho_2, \ldots, \rho_{M-1}, \xi_{M}+ \alpha- \sum_{j=1}^{M-1}(\eta_j-\xi_j) , \xi_{M+1}, \ldots >.$$
Then $\xi \le \rho$ by (\ref {e: 7}) and (\ref {e: 8}).
Moreover, by (\ref{e: 8}),
$\sum_{j=1}^n(\eta_j-\rho_j) \ge 0$ for all $n \le {M-1}$ and $\sum_{j=1}^{M-1}(\eta_j-\rho_j) = 0$, while for $n \ge M$,
\begin{align*}
\sum_{j=1}^n(\eta_j-\rho_j) &= \sum_{j=1}^{M-1}(\eta_j-\rho_j)+\eta_{M}-\Big(\xi_{M}+ \alpha -\sum_{j=1}^{M-1}(\eta_j-\xi_j)\Big) + \sum_{j=M+1}^n(\eta_j-\xi_j)\\
&= \sum_{j=1}^{n}(\eta_j-\xi_j)-\alpha >0
\end{align*}
by (\ref {e: 8}) and (\ref {e: 7}) and from this also $\varliminf_n\sum_{j=1}^n (\eta_j-\rho_j) = 0$.
Combining this with a check that $\rho$ is monotone nonincreasing completes the proof that $\rho \preccurlyeq \eta$.
\begin{alignat*}{3}
\rho_{M-1}& = \sum_{j=1}^{M-1}\rho_j - \sum_{j=1}^{{M-2}}\rho_j \ge \sum_{j=1}^{M-1}\eta_j - \sum_{j=1}^{{M-2}}\eta_j =\eta_{M-1} &\qquad \text{(by (\ref {e: 8}))}\\
&\ge \eta_{M}> \xi_{M}+ \alpha - \sum_{j=1}^{M-1}(\eta_j-\xi_j)= \rho_{M}&\qquad \text{(by (\ref {e: 7}))}\\
&\ge \xi_{M}&\qquad \text{(since $\xi \le \rho$)}\\
& \ge \xi_{M+1}= \rho_{M+1}
\end{alignat*}

\ep

The extension of (A)(ii) to infinite sequences is not based on the finite sequence result but requires a different approach:
we recursively reduce $\eta$ in a controlled way so that at each step the reduced sequence is also monotone nonincreasing and majorizes $\xi$.
For this we use the following device. For $\eta\in \co*$ and $t \ge 0$, define the $\co*$-sequence
\be{e: 9}
\eta(t):= ~<\min \{t, \eta_n\}>.
\ee
Notice that
\be {e: 11} \xi \prec \eta\quad\text{and}\quad \sum_{j=1}^n(\eta_j-\xi_j)=0 \quad \Rightarrow \quad \xi^{(n)}\prec \eta^{(n)}.\ee

\bL{L:3.5}

Let $\xi, \eta \in \co*$ with $\xi \prec \eta$.
There is a minimal value $0 \le t \in [\xi_1,\eta_1]$ for which $\xi\prec \eta(t)$. \linebreak
For that value $\inf_k \sum_{j=1}^k(\eta(t)_j-\xi_j)=0$ and then either
\item [(i)] $\xi\preccurlyeq \eta(t)$\\
or
\item [(ii)] $\sum_{j=1}^n(\eta_j(t)-\xi_j)=0$ for some integer $n\ge 1$.
\eL
\bp
If $\xi_1=0$ and hence $\xi=0$, it suffices to choose $t=0$ for which $\eta(0) = 0$ satisfies the requirement.
Assume henceforth that $\xi_1>0$.
By definition, $\eta(\eta_1)=\eta$, hence $\xi \prec \eta(\eta_1)$.
If $s<\xi_1$, then $\eta(s)_1 = s < \xi_1 $ and hence $\xi \not\prec \eta(s)$.
Furthermore if $s, t\in [\xi_1,\eta_1]$ and $s< t$, then $ \eta(s) \le \eta(t)$, hence if $\xi \prec \eta(s)$ then also $\xi \prec \eta(t)$.
Thus $\{s \in [\xi_1, \eta_1] \mid \xi \prec \eta(s)\}$ is an interval. It is closed because
$$
\{s \in [\xi_1, \eta_1] \mid \xi \prec \eta(s)\} = \bigcap_{n=1}^\infty \{s \in [\xi_1, \eta_1]~\Big | ~\sum_{j=1}^n (\eta_j(s)-\xi_j)\ge0\}
$$
is closed, since for every $j$, $\eta_j(s)= \min\{s, \eta_j\}$ is a continuous function of $s$.
Thus
$$
\{s \in [\xi_1, \eta_1] \mid \xi \prec \eta(s)\} =[t, \eta_1] \quad \text{for some } t \in [\xi_1,\eta_1],
$$
i.e., $t$ is the minimal value of $s$ for which $\xi \prec \eta(s)$.

Assume by contradiction that $\inf_k \sum_{j=1}^k(\eta(t)_j-\xi_j)= \gamma >0$.
Let $N$ be the unique integer for which $$\eta_{N+1}< t\le \eta_{N}$$
and set
$$\epsilon:= \min \{\frac{\gamma}{N}, t- \eta_{N+1} \}$$
so $\eta_{N+1} \le t-\epsilon < \eta_{N}$. Then $\eta(t-\epsilon)\le \eta(t)\le \eta$ and $\eta(t-\epsilon)_j=\eta(t)_j= \eta_j$ for $j>N$.
So one has for every $k$,
$$\sum_{j=1}^k \big(\eta(t-\epsilon)-\xi_j\big)= \sum_{j=1}^k \big(\eta(t)-\xi_j\big)-\begin{cases} k\epsilon &k\le N\\
N\epsilon &k> N\end{cases}\quad \ge \gamma - N\epsilon \ge 0.$$
Thus $\xi\prec \eta(t-\epsilon)$, against the minimality of $t$, which proves that $\inf_k \sum_{j=1}^k(\eta(t)_j-\xi_j)=0$.

To complete the proof we must show that if condition (i) of this lemma fails, i.e., $\varliminf_k\sum_{j=1}^k(\eta(t)_j-\xi_j) >0$, then condition (ii) must hold.
>From the definition of $\eta(t)$ and $N$ we have
\be{e: 12}
\sum_{j=1}^k(\eta(t)_j-\xi_j) =\begin{cases} \sum_{j=1}^k (t- \xi_j) &\quad \text{for } 1\le k \le N\\
\sum_{j=1}^k(\eta_j- \xi_j) - \sum_{j=1}^{N}(\eta_j- t) &\quad \text{for } k> N
\end{cases}
\ee
and so
$$\sum_{j=1}^{N}(\eta_j- t)<\varliminf_k\sum_{j=1}^k(\eta_j-\xi_j).$$
By the monotonicity of $\xi$, $\sum_{j=1}^k (t- \xi_j) \ge t-\xi_1$ for every $1\le k \le N$.
Thus
$$
0=\inf_k \sum_{j=1}^k(\eta(t)_j-\xi_j)= \min\Big\{t-\xi_1, \inf_{k>N}\sum_{j=1}^k(\eta_j- \xi_j) - \sum_{j=1}^{N}(\eta_j- t) \Big\}.
$$
If $t=\xi_1$, then (ii) is satisfied by $n=1$.
If $t>\xi_1$, then
$$\inf_{k>N}\sum_{j=1}^k(\eta_j- \xi_j)= \sum_{j=1}^{N}(\eta_j- t)< \varliminf_k\sum_{j=1}^k(\eta_j-\xi_j).$$
From this strict inequality between the $\inf$ and the $\varliminf$ it follows that $\{ \sum_{j=1}^k(\eta_j-\xi_j)\} _{ k > N}$ attains a minimum at some $n>N$,
in which case then $\sum_{j=1}^n(\eta_j- \xi_j)= \sum_{j=1}^{N}(\eta_j- t)$.
But then by (\ref {e: 12}) when $k = n$, $\sum_{j=1}^{n}(\eta(t)_j-\xi_j) =0$ which establishes condition (ii).
\ep

\bT{T:3.6}
If $\xi, \eta \in \co*$ and $\xi \prec \eta$, then there is some $\zeta \in \co*$ for which $\xi\preccurlyeq \zeta \le \eta$.
\eT
\bp
By Lemma \ref {L:3.5} there is a $\xi_1\le t_1\le \eta_1$ for which $\xi\prec \eta(t_1)$ and either $\xi\preccurlyeq \eta(t_1)$ or $\sum_{j=1}^{n_1}(\eta(t_1)_j-\xi_j)=0$ for some integer $n_1\ge 1$.
Set $\zeta(1):=\eta(t_1)$ and observe that by definition, $\zeta(1)\le \eta$.
In the first case, i.e., when $\xi\preccurlyeq \eta(t_1)$, choosing $\zeta:= \zeta(1)$ satisfies the claim.

In the second case, i.e., when $\sum_{j=1}^{n_1}(\zeta(1)_j-\xi_j)=0$ for some integer $n_1$ and $\xi\prec \eta(t_1) = \zeta(1)$,
it follows using (\ref{e: 11}) that $\xi^{(n_1)}\prec \zeta(1)^{(n_1)}$.
Apply again Lemma \ref {L:3.5} to these two truncated sequences and obtain a $t_2$ with
$$\xi_{n_1+1}=\xi^{(n_1)}_1\le t_2\le \zeta(1)^{(n_1)}_1 \le t_1$$ for which $\xi^{(n_1)}\prec \zeta(1)^{(n_1)}(t_2)$ and
either $\xi^{(n_1)}\preccurlyeq \zeta(1)^{(n_1)}(t_2)$ or $\sum_{j=1}^{n_2-n_1}\big( \zeta(1)^{(n_1)}(t_2)_j-\xi^{(n_1)}_j\big)=0$ for some integer $n_2> n_1$.
Observe that since $t_2\le t_1$, it follows that $\zeta(1)(t_2)=\eta(t_2)$ and $$\zeta(1)^{(n_1)}(t_2)=\zeta(1)(t_2)^{(n_1)}=\eta(t_2)^{(n_1)}.$$

Now splice together the beginning segment of the sequence $\zeta(1)$ with the sequence $ \zeta(1)^{(n_1)}(t_2)$ as follows:
$$
\zeta(2)_j:=\begin{cases}\zeta(1)_j = \eta(t_1)_j&1\le j\le n_1\\
\zeta(1)(t_2)_j =\eta(t_2)_j& j > n_1.
\end{cases}
$$
Thus $\zeta(2)$ is a monotone nonnegative sequence and $\zeta(2)\le \zeta(1)\le \eta$. Moreover, if $ 1\le n\le n_1$ then
$$\sum_{j=1}^n(\zeta(2)_j-\xi_j)=\sum_{j=1}^n(\eta(t_1)_j-\xi_j) \ge 0$$ and if $n>n_1$, then
\begin{align*}
\sum_{j=1}^n(\zeta(2)_j-\xi_j)&=\sum_{j=1}^{n_1}(\eta(t_1)_j-\xi_j)+ \sum_{j=n_1+1}^{n}(\eta(t_2)_j-\xi_j)\\
&= \sum_{j=n_1+1}^{n}(\eta(t_2)_j-\xi_j)\\
&= \sum_{j=1}^{n-n_1}(\eta(t_2)^{(n_1)}_j-\xi^{(n_1)}_j) \\
&= \sum_{j=1}^{n-n_1}(\zeta(1)^{(n_1)}_j(t_2)-\xi^{(n_1)}_j) \ge 0.
\end{align*}
Thus $\xi\prec \zeta(2)$. Furthermore, $\xi^{(n_1)}\preccurlyeq \zeta(1)^{(n_1)}(t_2)$ if and only if
\ba \varliminf_n\sum_{j=1}^n\big( \zeta(1)^{(n_1)}(t_2)-\xi^{(n_1)}_j \big) &= \varliminf_n\sum_{j=n_1+1}^{n_1+n}\big( \zeta(1)(t_2)_j-\xi_j \big) \\
&=\varliminf_n\sum_{j=n_1+1}^{n_1+n}\big( \eta(t_2)_j-\xi_j \big)\\
&=\varliminf_n\sum_{j=1}^n(\zeta(2)_j-\xi_j)=0,
\end{align*}
that is, if and only if $\xi\preccurlyeq \zeta(2)$.
In this case, we choose $\zeta:=\zeta(2)$. If however $\xi\not\preccurlyeq \zeta(2)$, then there must be an integer $n_2>n_1$ for which
$\sum_{j=1}^{n_2}(\zeta(2)_j-\xi_j)=0$ and hence by (\ref {e: 11}), $\xi^{(n_2)}\prec \zeta^{(n_2)}$ and we can repeat the construction.

There are now two possibilities.
In the first case, after a finite number of steps, say $M$ steps, we construct a nonnegative monotone sequence $\zeta(M)$ with $\xi\preccurlyeq \zeta(M)\le \eta$ thus satisfying the theorem.

In the second case, the process never terminates and so we obtain an infinite monotone collection of $\co*$-sequences
$\zeta(k)\ge \zeta(k+1)$ with $\xi\prec \zeta(k)$ and a strictly increasing sequence of integers $n_k$ with $n_o=0$ and $\zeta(0):=\zeta$ for which
\begin{itemize}
\item [(i)] $ \zeta(k+1)_j= \zeta(k)_j$ for every $k$ and every $1\le j\le n_k$;
\item [(ii)] $\sum_{1}^{n_k}\big(\zeta(k)_j-\xi_j\big)=0$ for every $k$.
\end{itemize}
Then splice together the sequences $\zeta(k)$ into a sequence $\zeta$ as follows:
$$\zeta_j:=\begin{cases}\zeta(1)_j &1\le j \le n_1\\
\zeta(2)_j& n_1< j\le n_2\\
\cdots &\phantom{aaa}\cdots
\end{cases}
$$
It is immediate to verify that $\zeta\in \co*$, $\zeta \le \eta$, $\xi\prec \zeta$ and furthermore, $\sum_{j=1}^{n_k}(\zeta_j-\xi_j)=0$ for every $k$. Thus \be{e: 13} \xi\prec_b \zeta \le \eta\ee and hence $\xi\preccurlyeq \zeta\le \eta$.
\ep

In fact condition (\ref {e: 13}) is equivalent to non-termination of the above construction, which is the second case just described.
Indeed, observe that if $\xi \prec \eta$ satisfy condition (ii) of Proposition \ref {P:3.2},
then the minimum of $\sum_{j=1}^k(\eta_j- \xi_j)$ in the proof of Lemma \ref {L:3.5} being attained, case (ii) of Lemma \ref {L:3.5} holds.
Thus in the above proof of Theorem \ref {T:3.6}, there is an $n_1$ for which $\sum_{j=1}^{n_1}(\zeta(1)_j-\xi_j)=0$.
Since for every $k\ge 1$
$$\sum_{j=1}^{k}(\zeta(1)^{(n_1)}_j-\xi^{(n_1)}_j) = \sum_{j=n_1+1}^{n_1+k}(\zeta(1)_j-\xi_j )= \sum_{j=1}^{n_1+k}(\zeta(1)_j-\xi_j ),$$
it follows that $\big\{\sum_{j=1}^{k}(\zeta(1)^{(n_1)}_j-\xi^{(n_1)}_j) \big\}_{k>m}$ also attains a minimum for every $m$,
i.e., the pair of sequences $\xi^{(n_1)}\prec \zeta(1)^{(n_1)}$ satisfies the same condition (ii) of Proposition \ref {P:3.2} as did the original pair $\xi,\eta$. Iterating this reasoning, the same holds at every step of the proof.
As a consequence, the construction in the proof of Theorem \ref {T:3.6} does not terminate after finitely many steps but yields an infinite strictly increasing sequence of integers $n_k$ for which $\sum_{j=1}^{n_k}(\zeta_j-\xi_j)=0$, i.e., the sequence $\zeta$ block-majorizes $\xi$ (see(\ref {e: 13})).
The converse is elementary (e.g., see the proof of implication (iii) $\Rightarrow$ (i) of Proposition \ref {P:3.2}). To summarize:

\bC{C:3.7}
Let $\xi, \eta \in \co*$ and $\xi \prec \eta$.
Then the condition $\xi \prec_b\zeta \le \eta$ for some $\zeta \in \co*$ can be added to the list of equivalent conditions of Proposition \ref {P:3.2}.
\eC

\bR{R:3.8}

Let $\xi, \eta \in \co*$ with $\xi\preccurlyeq \eta$.
Then $\xi \prec \zeta \le \eta$ with $\zeta \in \co*$ implies $\zeta = \eta$, i.e., there are no proper such $\zeta$.
Likewise $\xi \le \rho \prec \eta$ with $\rho \in \co*$ implies $\xi = \rho$.
Indeed, $$\eta_n-\zeta_n \le \sum _{j=1}^m(\eta_j-\zeta_j) \le \sum _{j=1}^m(\eta_j-\xi_j) \quad \text{for every }m\ge n,$$
and hence $\eta_n-\zeta_n \le {\varliminf}\sum_{j=1}^m ( \eta_j -\xi_j) =0$ for every $n$.
A similar argument shows that $\rho_n-\xi_n \le {\varliminf}\sum_{j=1}^m ( \eta_j -\xi_j) =0$ for every $n$.
\eR

Next we proceed with the extension to infinite sequences of Propositions \ref{P:3.1}(B)(i)-(ii).
In the special case when $\sum_{j=1}^\infty \xi_j = \sum_{j=1}^\infty\eta_j < \infty$, by (\ref{e: 2}) one has $\eta \prec_\infty\xi \Leftrightarrow \xi \prec \eta$ from which the results in Proposition \ref{P:3.2} and Theorems \ref {T: 3.4} and \ref {T:3.6} can be immediately reformulated in terms of majorization at infinity.
However, a more general extension of these results to both summable and nonsummable sequences requires another approach.

For this we use another device which in a sense is the mirror image of the one in (\ref{e: 9}).
For $\xi\in \co*$, $t\ge0$ and $p\in \mathbb N$, set
\be{e: 14}
\xi (t,p) :=
\begin{cases}
<t, t, \ldots> & \text {if $p=1$}\\
<\xi_1, \xi_2, \ldots, \xi_{p-1}, t, t, \ldots> & \text {if $p>1.$}
\end{cases}
\ee
Despite that the infinite sequence $\xi(t,p)$ is monotone nonincreasing only when $p=1$ or $t\le \xi_{p-1}$ and does not converge to $0$ unless $t=0$, we will in any case continue to use the majorization notation $\eta\prec \xi (t,p)$ to denote the property that
$\sum_{j=1}^n( \xi (t,p)_j-\eta_j)\ge 0$ for all $n$.

\bL{L:3.9}
Let $\xi, \eta\in \co*$ and assume that $\eta \not\prec \xi$. Then there are integers $n_1\ge p_1\ge1$ and a positive number $t_1$ with
$\begin{cases}t_1=\eta_1,~ p_1=1 & \text{if }\eta_1> \xi_1\\
\xi_{p_1} < t_1\le \xi_{p_1-1}, ~p_1>1 & \text{if } \eta_1\le \xi_1
\end{cases}$ \quad for which
\item[(i)] $\xi(t_1,p_1)$ is monotone nonincreasing;
\item [(ii)] $\xi\le \xi(t_1,p_1)$;
\item [(iii)] $\eta \prec \xi(t_1,p_1)$;
\item [(iv)] $\sum_{j=1}^{n_1}(\xi (t_1, p_1)_j-\eta_j)=0$;
\item [(v)] $\xi(t_1, p_1)_{1}= \max\{ \xi_{1}, \eta_{1} \}$
\item [(vi)] $\xi(t_1, p_1)_{n_1}\ge \max\{ \xi_{n_1+1}, \eta_{n_1+1} \}.$
\eL
\bp
If $\eta_1> \xi_1$, set $n_1=p_1=1$ and $t_1=\eta_1$.
Then $\xi(t_1, p_1)=~<\eta_1,\eta_1, \cdots >$ and conditions (i)-(vi) are trivially satisfied.
So henceforth assume that $\eta_1\le \xi_1$.
Since $\eta \not \prec \xi$, it follows that for some $1 < N \in \mathbb N$, one has 
$$\sum_{j=1}^N\eta_j > \sum_{j=1}^N\xi_j = \sum_{j=1}^N\xi(\xi_N,N)_j,$$
i.e., $\eta \not\prec \xi(\xi_N,N)$, so the set $\{ p \in \mathbb N\mid \eta \not\prec \xi (\xi_p,p)\}$ is nonempty. Set
\be{e: 13'}
p_1:= \min \{ p \in \mathbb N\mid \eta \not\prec \xi (\xi_p,p)\}.
\ee
Equivalently, $p_1$ is the minimum $n$ for which $\sum_{j=1}^n\eta_j > \sum_{j=1}^n\xi_j$.
Then $p_1\le N$ and since $\eta_1\le \xi_1$ and hence $\eta \prec \xi (\xi_1,1)= ~<\xi_1, \xi_1, \cdots>$, it follow that $p_1>1$.
Thus for arbitrary $t \ge 0$, $\xi(t, p_1)_{1}=\xi_1= \max\{ \xi_{1}, \eta_{1} \}$ and condition (v) is satisfied.
For every $t\ge 0$,
\be{e: 15}
\sum_{j=1}^n(\xi
(t,p_1)-\eta_j)=
\begin{cases}
\sum_{j=1}^n(\xi_j-\eta_j)\quad &\text{for } 1 \le n < p_1\\
\sum_{j=1}^{p_1-1}\xi_j + (n-p_1+1)t - \sum_{j=1}^n\eta_j &\text{for } n \ge p_1.
\end{cases}
\ee
Since $p_1 > 1$ and $\eta \prec \xi (\xi_{p_1-1}, p_1-1)$ by the minimality of $p_1$, for all $1 \le n < p_1$ one has
$$\sum_{j=1}^n(\xi_j-\eta_j)= \sum_{j=1}^n(\xi
(\xi_{p_1-1}, p_1-1)_j-\eta_j)\ge 0.$$
Thus
\be{e: 16}\sum_{j=1}^n(\xi
(t,p_1)-\eta_j)\ge 0\quad \text{for all $1 \le n < p_1$ and all $t\ge 0.$}
\ee
To identify $t_1$ we use the function defined for $n\ge p_1$:
$$
t(n):= \frac{1}{n-p_1+1}\big( \sum_{j=1}^n\eta_j -\sum_{j=1}^{p_1-1}\xi_j \big)= \frac{n}{n-p_1+1}(\eta_a)_n - \frac{1}{n-p_1+1}\sum_{j=1}^{p_1-1}\xi_j
$$
where $\eta_a := ~<(\eta_a)_n>$ is the arithmetic mean sequence of $\eta$.
It is well know that if $\eta\in \co*$ then also $\eta_a\in \co*$, i.e., $(\eta_a)_n\downarrow 0$ , hence $t(n)\to 0$.
By definition of $p_1$, $\eta \not\prec \xi(\xi_{p_1},p_1)$ and $\sum_{j=1}^{p_1}\eta_j > \sum_{j=1}^{p_1}\xi_j\ge \sum_{j=1}^{p_1-1}\xi_j$,
so it follows that $t(p_1) >0$.
Since $t(n) \rightarrow 0$ and $t(p_1) >0$, one has for some $n_1 \ge p_1$, $$t_1:=t(n_1) = \max_{n\ge p_1} t(n) > 0.$$
Thus
\be {e: 17} (n-p_1+1)t_1 \ge \sum_{j=1}^n\eta_j - \sum_{j=1}^{p_1-1}\xi_j \quad \text{for all } n \ge p_1, \text{ and with equality for }n=n_1.
\ee
But then from (\ref {e: 15}), (\ref {e: 16}), and (\ref{e: 17}) we have that $\sum_{j=1}^n(\xi(t_1,p_1)-\eta_j)\ge 0$ for all $n\ge 1$,
which establishes (iii), and that
$
\sum_{j=1}^{n_1}(\xi
(t_1,p_1)-\eta_j)= 0$, which establishes (iv).

Since $n_1 \ge p_1$, if $s< t_1$, then $\sum_{j=1}^{n_1}(\xi (s,p_1)_j-\eta_j)< 0$ hence $\eta\not\prec \xi (s,p_1)$.
But since $$\eta \prec \xi (\xi_{p_1-1}, p_1-1) = \xi (\xi_{p_1-1}, p_1),$$ it follows then that $\xi_{p_1-1} \ge t_1$.
On the other hand, from $\eta \not \prec \xi(\xi_{p_1}, p_1)$ and by (iii): $\eta \prec \xi(t_1, p_1) \prec \xi(t, p_1)$ for all $t \ge t_1$,
it follows that $\xi_{p_1} < t_1$.
So $\xi_{p_1} < t_1 \le \xi_{p_1-1}$.
This implies by Definition (\ref {e: 14}) both that $\xi (t_1, p_1)$ is monotone nonincreasing and that $\xi\le \xi (t_1, p_1)$, thus establishing (i) and (ii).

Furthermore, from $p_1\le n_1$ it follows that $\xi (t_1, p_1)_{n_1}=t_1 > \xi_{p_1}\ge \xi_{n_1}\ge \xi_{n_1+1}$ and from (iii): $\eta \prec \xi(t_1,p_1)$ and (iv) it follows that
\[
\xi
(t_1, p_1)_{n_1}\ge \xi
(t_1, p_1)_{n_1+1}= \sum _{j=1}^{n_1+1}\xi
(t_1, p_1)_j- \sum _{j=1}^{n_1}\xi
(t_1, p_1)_j
\ge \sum _{j=1}^{n_1+1}\eta_j - \sum _{j=1}^{n_1}\eta_j = \eta_{n_1+1}
\]
establishing (vi) and concluding the proof.
\ep
\bL{L:3.10.0}
Let $\xi, \eta \in \co*$ and assume that $\eta^{(n)} \not \prec \xi^{(n)}$ for $n=0, 1, 2\ldots$ Then $\zeta \le \eta$ and $\zeta\prec_b \xi$ for some $\zeta \in \co*$.
\eL
\bp
Since $\eta \not \prec \xi $, there is some $n$ for which $\sum_{j=1}^n(\eta_j -\xi_j) >0$.
Let $n_1$ be the smallest such integer.
If $n_1=1$, then the inequality $\sum_{j=n}^{n_1}(\eta_j- \xi_j )\ge 0$  holds trivially (and strictly) for $0<n\le n_1$.
If $n_1>1$, $\sum_{j=1}^{n}(\eta_j -\xi_j)\le 0$ for every $0< n< n_1$.
But then for $2\le n\le n_1$ we have
$$\sum_{j=n}^{n_1}(\eta_j -\xi_j) = \sum_{j=1}^{n_1}(\eta_j -\xi_j) - \sum_{j=1}^{n-1}(\eta_j -\xi_j) >0. $$
The same inequality holds by definition of $n_1$ for $n=1$.
Since by hypothesis $\eta^{(n_1)} \not\prec \xi^{(n_1)}$, we can iterate this construction and obtain a strictly increasing sequence of integers $ n_k$ starting with $n_o=0$ for which
$\sum_{j=n}^{n_k}(\eta_j- \xi_j )\ge 0$ for all $n_{k-1} < n \le n_k$ and all $k\in \mathbb N$.

But then
\[
<\xi_{n_{k-1}+1}, \xi_{n_{k-1}+2},\ldots, \xi_{n_k}> \: \prec_\infty\: <\eta_{n_{k-1}+1}, \eta_{n_{k-1}+2},\ldots, \eta_{n_k}>.
\]
Thus by Proposition \ref {P:3.1}(B)(ii) and (\ref {e: 2}) there is for each $k$ a monotone nonincreasing finite sequence
$$<\zeta_{n_{k-1}+1}, \zeta_{n_{k-1}+2},\ldots, \zeta_{n_k}> \: \le \: <\eta_{n_{k-1}+1}, \eta_{n_{k-1}+2},\ldots, \eta_{n_k}> $$
with
$$<\zeta_{n_{k-1}+1}, \zeta_{n_{k-1}+2},\ldots, \zeta_{n_k}>\: \preccurlyeq \: <\xi_{n_{k-1}+1}, \xi_{n_{k-1}+2},\ldots, \xi_{n_k}>$$
and so in particular $\xi_{n_k}\le \zeta_{n_k}$.
Let $\zeta$ be the sequence obtained by splicing together the finite sequences thus obtained for each interval $n_{k-1} < n \le n_k$.
Then $\zeta\le \eta$ and if $ n_{k-1} < n \le n_k$ for $k>1$, then
\[
\sum_{j=1}^{n} \zeta_j = \sum_{i=1}^{k-1}\Big(\sum_{j= n_{i-1}+1}^{n_i}\zeta_j \Big) + \sum_{j= n_{k-1}+1}^n\zeta_j \le \sum_{i=1}^{k-1}\Big(\sum_{j=n_{i-1}+1}^{n_i}\xi_j \Big) + \sum_{j= n_{k-1}+1}^n\xi_j = \sum_{j=1}^{n} \xi_j
\]
with equality holding for $n = n_k$. The same conclusion holds directly for $k=1$. Thus, $\zeta\prec_b \xi$. Finally,
\[
\zeta_{n_k +1} = \sum_{j=1}^{n_k +1}\zeta_j- \sum_{j=1}^{n_k }\zeta_j\le \sum_{j=1}^{n_k +1}\xi_j- \sum_{j=1}^{n_k }\xi_j = \xi_{n_k +1}\le \xi_{n_k}\le \zeta_{n_k}
\]
which proves that $\zeta$ is monotone nonincreasing.
\ep
This lemma and the ideas in its proof lead to the following characterization of the case when there is an intermediate sequence $\zeta$ for which $\zeta \le \eta$ and $\zeta\prec_b \xi$.

\bP{P:3.10} Let $\xi, \eta \in \co*$. Then the following conditions are equivalent.
\item[(i)] There is a strictly increasing sequence of integers $ n_k$ starting with $n_o=0$ for which
$\sum_{j=n}^{n_k}(\eta_j- \xi_j )\ge 0$ for all $n_{k-1} < n \le n_k$ and all $k\in \mathbb N$.
\item[(ii)] Either
\item[(a)] $\eta^{(n)} \not \prec \xi^{(n)}$ for $n=0, 1, 2\ldots$\\
or
\item[(b)] There is an $N\in \mathbb N$ for which $\eta^{(N)} \prec_b \xi^{(N)}$ and $\sum_{j=n}^N(\eta_j-\xi_j) \ge 0$ for all $1 \le n \le N$.
\item [(iii)] $\zeta \le \eta$ and $\zeta\prec_b \xi$ for some $\zeta \in \co*$.
\item [(iii$'$)] $\eta\prec_b \rho$ and $\xi\le \rho$ for some $\rho \in \co*$ .
\eP
\bp
First we prove that (i) and (iii) are equivalent. 
\item(i) $\Rightarrow$ (iii)
This implication is the second part of the proof of Lemma \ref {L:3.10.0}.
\item(iii) $\Rightarrow$ (i) Let $n_k$ be the strictly increasing sequence of indices for which $\sum_{j=1}^{n_k}(\xi_j-\zeta_j)=0$
or, equivalently, $\sum_{j=n_{k-1}+1}^{n_k}(\xi_j-\zeta_j)=0$ for all $k$.
For each $n \ge 1$, choose $k$ for which $n_{k-1}<n\le n_k$.
Then since $\zeta \prec \xi$, for $n_{k-1}+1<n\le n_k$ one has
$$
\sum_{j=n}^{n_k}(\zeta_j-\xi_j)= \sum_{j=n_{k-1}+1}^{n_k}(\zeta_j-\xi_j)- \sum_{j=n_{k-1}+1}^{n-1}(\zeta_j-\xi_j)= \sum_{j=n_{k-1}+1}^{n-1}( \xi_j -\zeta_j)= \sum_{j=1}^{n-1}( \xi_j -\zeta_j)\ge 0,
$$
and for $n=n_{k-1}+1$, $\sum_{j=n}^{n_k}(\zeta_j-\xi_j)= 0$ by hypothesis, that is, $\sum_{j=n}^{n_k}(\zeta_j-\xi_j) \ge 0$ for $n_{k-1}<n\le n_k$.
Thus since $\zeta \le \eta$ by hypothesis, for all $n_{k-1}<n\le n_k$,
$$
\sum_{j=n}^{n_k}(\eta_j-\xi_j)\ge \sum_{j=n}^{n_k}(\zeta_j-\xi_j)\ge 0.
$$
Next we prove the equivalence of (i), (ii), and (iii$'$).
\item (i) $\Rightarrow$ (ii)
Assuming (i) holds, then
\be{e: 17'}
\sum _{j=n}^{n_k}(\eta_j-\xi_j) \ge 0 ~\text{for every integer}~ n \le n_k~ \text{and every}~ k.
\ee
If there are infinitely many indices $k$ for which $\sum _{j=n_{k}+1}^{n_{k+1}}(\eta_j-\xi_j) >0$,
then for every integer $n\ge 0$ there is an index $k$ for which $n_{k-1} > n$ and $\sum_{j=n_{k-1}+1}^{n_k}(\eta_j-\xi_j) > 0$ so
$$\sum_{j=1}^{n_k-n}(\eta_j^{(n)}-\xi_j^{(n)}) ~= \sum_{j=n+1}^{n_k}(\eta_j-\xi_j) ~= \sum_{j=n+1}^{n_{k-1}}(\eta_j-\xi_j) ~+ \sum_{j=n_{k-1}+1}^{n_k}(\eta_j-\xi_j) ~>~ 0.$$
Thus $\eta^{(n)} \not \prec \xi^{(n)}$ and hence (ii)(a) holds.

If on the other hand $\sum _{j=n_{k}+1}^{n_{k+1}}(\eta_j-\xi_j) =0$ for all $k\ge K$ for some $K\in \mathbb N$, set $N=n_K$.
By (\ref{e: 17'}), $\sum _{j=n}^N(\eta_j-\xi_j) \ge 0$ for all integers $n \le N$.
Then for every $n> N$, letting $h\ge K$ be the index for which $n_h\le n< n_{h+1}$, one has
$$
\sum _{j=N+1}^n(\eta_j-\xi_j) =\sum _{j=N+1}^ {n_{h+1}}(\eta_j-\xi_j) - \sum _{j=n+1}^{ n_{h+1}}(\eta_j-\xi_j) = - \sum _{j=n+1}^ {n_{h+1}}(\eta_j-\xi_j) \le 0.
$$
Thus $\eta^{(N)} \prec \xi^{(N)}$. Furthermore, for every $k>K$,
$$
\sum_{j=1}^{n_k-N} (\xi^{(N)}_j- \eta^{(N)}_j) =0,
$$
hence $\eta^{(N)} \prec_b \xi^{(N)}$, and thus (ii)(b) holds.
\item (ii)(a) $\Rightarrow$ (iii$'$)
Since $\eta=\eta^{(0)}\not\prec \xi^{(0)}=\xi$, by Lemma \ref {L:3.9},
there is a sequence $\xi(t_1,p_1)$ and an integer $n_1\ge p_1$ satisfying conditions (i)-(vi) of Lemma \ref {L:3.9}.
Define $$\rho_j:= \xi(t_1,p_1)_j\quad\text{for}\quad 1\le j\le n_1.$$
Thus in particular $<\rho_1, \rho_2, \cdots,\rho_{n_1}> $ is a monotone nonincreasing finite sequence for which
\ba
<\xi_1, \xi_2, \cdots, \xi_{n_1}>&~\le~ <\rho_1, \rho_2, \cdots,\rho_{n_1}>\\
<\eta_1, \eta_2, \cdots, \eta_{n_1}>&~\preccurlyeq ~<\rho_1, \rho_2, \cdots,\rho_{n_1}>\\
\rho_1&=\max\{\xi_1, \eta_1\}\\
\rho_{n_1}&~\ge~ \max\{ \xi_{n_1+1}, \eta_{n_1+1} \}.
\end{align*}
By hypothesis, $\eta^{(n_1)}\not\prec \xi^{(n_1)}$.
Thus applying again Lemma \ref {L:3.9} we find a sequence $\xi^{(n_1)}(t_2,p_2)$ and an integer $n_2 \ge n_1 + p_2$ satisfying conditions (i)-(vi) in Lemma \ref {L:3.9}. Explicitly,
$$
\sum _{j=1}^{n_2-n_1}\big(\xi^{(n_1)}(t_2,p_2)_j-\eta^{(n_1)}_j\big)=0
$$
and
$$
\xi^{(n_1)}(t_2,p_2)_{n_2-n_1}\ge \max\{ \xi^{(n_1)}_{n_2-n_1+1}, \eta^{(n_1)}_{n_2-n_1+1} \}= \max\{ \xi_{n_2+1}, \eta_{n_2+1} \}.
$$
Thus defining
$$\rho_j:= \xi^{(n_1)}(t_2,p_2)_{j-n_1}\quad\text{for}\quad n_1< j\le n_2$$
we obtain again
a monotone nonincreasing finite sequence for which
\ba
<\xi_{n_1+1}, \xi_{n_1+2}, \cdots, \xi_{n_2}>&~\le~ <\rho_{n_1+1}, \rho_{n_1+2}, \cdots,\rho_{n_2}>\\
<\eta_{n_1+1}, \eta_{n_1+2}, \cdots, \eta_{n_2}>&~\preccurlyeq ~<\rho_{n_1+1}, \rho_{n_1+2}, \cdots,\rho_{n_2}>\\
\rho_{n_1+1}&=\max\{\xi_{n_1+1}, \eta_{n_1+1}\}\\
\rho_{n_2}&~\ge~ \max\{ \xi_{n_2+1}, \eta_{n_2+1} \}.
\end{align*}
Iterating this construction obtains a strictly increasing sequence of indices $n_k$ and for each $k \in \mathbb N $ a monotone nonincreasing finite sequence
\ba
<\xi_{n_{k-1}+1}, \xi_{n_{k-1}+2}, \cdots, \xi_{n_{k}}>&~\le~ <\rho_{n_{k-1}+1}, \rho_{n_{k-1}+2}, \cdots,\rho_{n_{k}}>\\
<\eta_{n_{k-1}+1}, \eta_{n_{k-1}+2}, \cdots, \eta_{n_{k}}>&~\preccurlyeq ~<\rho_{n_{k-1}+1}, \rho_{n_{k-1}+2}, \cdots,\rho_{n_{k}}>\\
\rho_{n_k}&~\ge~ \max\{ \xi_{n_k+1}, \eta_{n_k+1} \}\\
\rho_{n_{k}+1}&=\max\{\xi_{n_{k}+1}, \eta_{n_k+1}\}.
\end{align*}
Thus $\xi\le \rho$ and $\eta \prec_b \rho$. Since $\rho_{n_k} \ge \rho_{n_k+1}$ it follows that $\rho$ is monotone nonincreasing and since $ \rho_{n_{k}+1}\to0$, it follows that $\rho\in \co*$.
\item (ii)(b) $\Rightarrow$ (iii$'$)
Let $\rho:=\,< \xi_1+\sum_{j=1}^N(\eta_j-\xi_j), \xi_2, \xi_3, \ldots>$. Then $\xi \le \rho \in \co*$ and
\[
\sum_{j=1}^n(\rho_j - \eta_j)= \sum_{j=1}^n(\xi_j - \eta_j) +\sum_{j=1}^N(\eta_j-\xi_j)=
\begin{cases}
\sum_{j=n+1}^N(\eta_j-\xi_j)\ge 0 \quad &\text{for } 1\le n < N\\
0&\text{for } n = N\\
\sum_{j=N+1}^n(\xi_j - \eta_j)\ge 0 &\text{for } n > N.
\end{cases}
\]
Thus $\eta \prec \rho$ and it is also immediate to see that $\eta \prec_b \rho$.
 \item (iii$'$) $\Rightarrow$ (i) Similar proof as in (iii) $\Rightarrow$ (i).
\ep

\bT{T:3.11} Let $\xi,\eta \in (\ell^1)^*$ and $ \xi \prec_\infty \eta$.
\item [(i)] $\xi \preccurlyeq_\infty \zeta \le \eta$ for some $\zeta \in \co*$.
\item [(ii)] $\xi \le \rho\preccurlyeq_\infty \eta$ for some $\rho \in \co*$.
\eT

\bp
\item [(i)] If $\eta^{(n)} \not \prec \xi^{(n)}$ for all $n$, then by Lemma \ref {L:3.10.0} there is a $\zeta \in \co*$ for which $\zeta \le \eta$ and $\zeta\prec_b \xi$, hence $\zeta\preccurlyeq \xi$, and thus $\xi \preccurlyeq_\infty \zeta $.
Thus it remains to consider the case when $\eta^{(N)} \prec \xi^{(N)}$ for some $N$, i.e., $\sum_{j=N+1}^n(\xi_j-\eta_j) \ge 0$ for every $n \ge N+1$. Thus $\sum_{j=N+1}^\infty(\xi_j-\eta_j) \ge 0$. By hypothesis, $\sum_{j=N+1}^\infty(\xi_j-\eta_j) \le 0$, hence we have equality, i.e., $\eta^{(N)} \preccurlyeq \xi^{(N)}$, or, equivalently, $\xi^{(N)}\preccurlyeq_\infty \eta^{(N)}$. In particular, $\xi_{N+1}\ge \eta_{N+1}$. Moreover,
$$
\sum_{j=n}^N(\eta_j-\xi_j)= \sum_{j=n}^\infty(\eta_j-\xi_j)\ge 0\quad\text{for every }1\le n\le N,$$
i.e., $<\xi_1,\xi_2, \ldots, \xi_N> \prec_\infty <\eta_1,\eta_2, \ldots, \eta_N>$. By Proposition \ref {P:3.1}(B)(ii), there is a finite nonincreasing sequence $\zeta \in \mathbb (R^N)^+$ with $\xi \preccurlyeq_\infty \zeta \le \eta$ on the integer interval $1, 2, \ldots, N$. In particular, $\xi_N \le \zeta_N$. Define the infinite sequence $\zeta:=\,<\zeta_1, \zeta_2, \ldots, \zeta_N, \eta_{N+1}, \eta_{N+2}, \ldots>$.
Since $\zeta_{N} \ge \xi _{N }\ge \xi _{N +1}\ge \eta_{N +1}= \zeta_{N +1}$, it follows that $\zeta$ is monotone nonincreasing. It is now immediate to see that $\xi \preccurlyeq_\infty \zeta \le \eta$.\\
\item [(ii)] It suffices to choose $\rho:=\,<\xi_1+\sum_{j=1}^\infty(\eta_j-\xi_j), \xi_2, \xi_3, \ldots >$.
\ep

\section{Characterization of am-closed ideals} \label{S: 4}

In this section we consider operator ideals in a separable infinite-dimensional Hilbert space $H$. We fix an orthonormal basis and we consider the maximal abelian algebra $\mathscr D$ of $B(H)$ of the diagonal operators relative to this basis. $E$ denotes the conditional expectation onto $\mathscr D$, that is, the operation of ``taking the main diagonal".

To characterize the arithmetic mean closure $I^-:=~_a(I_a)$ of an operator ideal $I$, it suffices to characterize its positive diagonal operators (with respect to a fixed basis), or merely its positive non-increasing diagonal operators, that is, the sequences comprising its characteristic set $\Sigma(I^-)$.

Notice that a sequence $\xi\in \co*$ belongs to $\Sigma(I^-) =\Sigma(_a(I_a)) $ if and only if $\xi_a\le \eta_a$ for some $\eta\in \Sigma(I)$ (by definition) if and only if $\xi\prec \eta$ for some $\eta\in \Sigma(I)$ (by (\ref{e: 5})).
This proves the (elementary) equivalence of (i) and (ii) in the following proposition which provides also characterizations of $\Sigma(I^-)$ in terms of stochastic matrices.

\bP{P:4.1}
Let $I$ be an ideal and let $ \xi \in \co*$. Then the following conditions are equivalent.
\item[(i)] $\xi \in \Sigma (I^-)$.
\item[(ii)] $\xi \prec \eta$ for some $\eta \in \Sigma(I)$.
\item[(ii$'$)] $\xi \preccurlyeq \eta$ for some $\eta \in \Sigma(I)$.
\item[(iii)] $\xi =P \eta$ for some $\eta \in \Sigma(I)$ and some substochastic matrix $P$.
\item[(iii$'$)] $\xi =Q \eta$ for some $\eta \in \Sigma(I)$ and some orthostochastic matrix $Q$.
\item[(iv)] $\diag \xi = E(L \diag \eta\1 L^*)$ for some $\eta \in \Sigma(I)$ and some contraction $L$.
\item[(iv$'$)] $\diag \xi = E(U \diag \eta\1 U^*)$ for some $\eta \in \Sigma(I)$ and some orthogonal matrix $U$.\\
\noindent
If $I\supset \mathscr L_1$, then these conditions are all also equivalent to:
\item[(ii$''$)] $\xi\prec_b \eta$ for some $\eta \in \Sigma(I)$.
\item[(iii$''$)] $\xi =Q \eta$ for some $\eta \in \Sigma(I)$ and some block orthostochastic matrix $Q$.
\item[(iv$''$)] $\diag \xi = E(U \diag \eta\1 U^*)$ for some $\eta \in \Sigma(I)$ and some matrix $U$ a direct sum of finite orthogonal matrices.
\eP
\bp
\item[(i)] $ \Leftrightarrow$ (ii) Proven in the paragraph preceding this proposition.
\item[(ii)] $ \Rightarrow$ (ii$'$) Let $\xi \prec \eta$ for some $\eta \in \Sigma(I)$. By Theorem \ref{T:3.6}
there is a $\zeta \in \co*$ for which $\xi\preccurlyeq \zeta \le \eta$. By definition, $\Sigma(I)$ is hereditary, hence $\zeta \in \Sigma(I)$.
\item[(ii$'$)] $ \Rightarrow$ (iii$'$) by Theorem \ref {T:SH}.
\item[(iii$'$)] $ \Rightarrow$ (iv$'$) by Lemma \ref {L:2.3}.
\item [(iv$'$)]$ \Rightarrow$ (iv) Obvious.
\item [(iv)] $ \Rightarrow$ (iii) by Lemma \ref {L:2.3}.
\item [(iii)] $ \Rightarrow$ (ii) by \cite [Lemma 3.1] {aM64} (see Lemma \ref {L:Mark} herein).

\noindent Assume now that $I\supset \mathscr L_1$.
\item [(ii)] $\Rightarrow$ (ii$''$) If $\xi \in(\ell^1)^*$, then $\xi \in \Sigma (I)$ so there is nothing to prove.
If $\xi \not \in (\ell^1)^*$, then $\sum_{j=1}^\infty (2\eta_j-\xi_j)=\infty$. \linebreak
By Remark \ref {R:3.3}, condition (ii) of Proposition \ref{P:3.2} is satisfied for $\xi \prec 2\eta\in \Sigma(I)$
and by equivalent condition (iii) in Proposition \ref{P:3.2} there is a $\zeta \in \co*$ for which $\xi\prec_b \zeta \le 2\eta$ and so $\zeta \in \Sigma (I)$ by the hereditariness of $\Sigma (I)$.
\item [(ii$''$)] $\Rightarrow$ (ii) Obvious.
\item[(ii$''$)] $ \Leftrightarrow$ (iii$''$) by (\ref {e: 3}).
\item[(iii$''$)] $ \Leftrightarrow$ (iv$''$) by Lemma \ref {L:2.3}.
\ep

\bR{R:4.2} In order to prove that (i) $ \Rightarrow$ (ii$''$) we indeed need to assume that $I \supset \mathscr L_1$. Consider for instance the finite rank ideal $F$, which obviously does not contain $\mathscr L_1$. Then $ F \subsetneq F^-= \mathscr L_1$
by (\ref {e: 4}). On the other hand, it is clear that if a sequence is block majorized by a finitely supported sequence it must also be finitely supported. Thus no $\xi \in \Sigma (F^-)\setminus \Sigma (F)$ can satisfy condition (ii$''$).
\eR

The analogue of Proposition \ref {P:4.1} for am-$\infty$ closure is:

\bP{P:4.3}
Let $I\subset \mathscr L_1$ be an ideal and let $ \xi \in (\ell^1)^*$. Then the following conditions are equivalent.
\item[(i)] $\xi \in \Sigma (I^{-\infty}).$
\item[(ii)] $ \xi \prec_\infty\eta$ for some $\eta \in \Sigma(I).$
\item[(ii$'$)] $\xi \preccurlyeq_\infty \eta$ for some $\eta \in \Sigma(I).$
\item[(ii$''$)] $\eta \preccurlyeq \xi$ for some $\eta \in \Sigma(I).$
\item[(iii)] $P\xi \in \Sigma(I)$ for some column-stochastic matrix $P.$
\item[(iii$'$)] $Q\xi \in \Sigma(I)$ for some block orthostochastic matrix $Q.$
\item[(iv)] $E(V \diag \xi\1 V^*)\in I$ for some isometry $V.$
\item[(iv$'$)] $E(U \diag \xi\1 U^*)\in I$ for some matrix $U$ a direct sum of finite orthogonal matrices.
\eP
\bp
\item[(i)] $ \Leftrightarrow$ (ii) As in the proof of Proposition \ref {P:4.1}, this equivalence is a reformulation of the definition of $I^{-\infty}$.
\item[(ii)] $ \Leftrightarrow$ (ii$'$) One direction is by Theorem \ref{T:3.11}
and the same argument as in the proof of Proposition \ref {P:4.1}; the other direction is obvious.
\item[(ii$'$)] $ \Leftrightarrow$ (ii$''$) By (\ref {e: 2}).
\item[(ii$'$)] $ \Rightarrow$ (iii$'$) If $\eta_n>0$ for all $n$, then $\sum_{j=n}^\infty(2\eta_j-\xi_j)\ge \sum_{j=n}^\infty\eta_j >0$, hence $2\eta^{(n)} \not\prec \xi^{(n)}$ for every $n$ and thus condition (ii)(a) in Lemma \ref {L:3.10.0} applies to the pair of sequences $\xi, 2\eta$. Thus there is a $\zeta \in \co*$ for which $\zeta\prec_b \xi$ and $\zeta \le 2\eta$ and hence $\zeta \in \Sigma(I)$. By (\ref{e: 3}), $\zeta = Q\xi$ for some block orthostochastic matrix $Q$, i.e., $Q\xi\in \Sigma(I)$. If on the other hand $\eta_{N}=0$ for some $N\in \mathbb N$ and hence for all $n\ge N$, then also $\xi_{n}=0$ for all $n\ge N$ and hence by the Schur-Horn Theorem there is an $N\times N$ orthostochastic matrix $Q_0$ mapping $<\xi_1, \ldots, \xi_N>$ onto $<\eta_1, \ldots, \eta_N>$. But then, $Q:=Q_o\oplus I$ is also block orthostochastic and $Q\xi=\eta\in \Sigma(I)$.
\item[(iii$'$)] $ \Rightarrow$ (iv$'$) By Lemma \ref {L:2.3}.
\item[(iv$'$)] $ \Rightarrow$ (iv) Obvious.
\item[(iv)] $ \Rightarrow$ (iii) By Lemma \ref{L:2.3}.
\item[(iii)] $ \Rightarrow$ (ii$''$) By Lemma \ref {L:2.10}.
\ep

\bC{C:4.4}
For every ideal $I$, $ I^-\cap \mathscr D=E(I)$.
\eC
\bp
Since ideals are the linear span of their positive parts and the conditional expectation $E$ is linear, it is enough to prove that $(I^-)^+\cap \mathscr D=E(I^+) $.

Let $A\in I^+$. Let $s(A)$ denote the sequence of singular numbers of $A$, namely the sequence of eigenvalues of $A$ repeated by multiplicity and in nonincreasing order. By \cite [Proposition 6.4]{vKgW10-Maj} , $s(E(A)) \prec s(A)$ and since $s(A)\in \Sigma(I)$ by definition, by Proposition \ref{P:4.1} it follows that
$s(E(A)) \in \Sigma (I^-)$, i.e., $E(A)\in I^-$. This proves that $E(I) \subset I^- \cap \mathscr D$.

Let now $B\in I^- \cap \mathscr D$. Then $s(B)\in \Sigma (I^-)$, hence by Proposition \ref{P:4.1}, $\diag s(B)= E(U\diag \eta U^*)$ for some orthogonal matrix $U$ and some $\eta\in \Sigma (I)$. Now $B= \Pi\diag s(B)\Pi^*$ for some permutation matrix $\Pi$ and since $\Pi$ commutes with the conditional expectation $E$ we have:
$$ B=\Pi E(U\diag \eta U^*) \Pi^*= E(\Pi U\diag \eta U^*\Pi^*).$$
Now $\diag \eta \in I$, hence $\Pi U\diag \eta U^*\Pi^*\in I$ and thus
$B\in E(I)$, which concludes the proof.
\ep

\bT{T:4.5}
An ideal $I$ is am-closed if and only if it is diagonally invariant, that is, $E(I) \subset I$.
\eT
\bp
If $I$ is am-closed, i.e., $I=I^-$, then by Corollary \ref{C:4.4}, $E(I) = I\cap \mathscr D\subset I$.
Conversely, if $E(I) \subset I$ then, again by Corollary \ref {C:4.4}, $I^-\cap \mathscr D\subset I$.
But since every selfadjoint element of $I^-$ is diagonalizable, and since $I$ and $I^-$ are ideals and hence unitarily invariant,
this implies that $I^- \subset I$ and hence $I^- = I$, i.e., $I$ is am-closed.
\ep

Proposition \ref {P:4.1} leads naturally to the following notion of invariance of an ideal under the action of a class of substochastic matrices. Recall from the remarks after Definition \ref {D:stoch} that if $\xi\in\co*$ and $P$ is a substochastic matrix, then $P\xi \in \text {c}_\text{o}$ and hence $(P\xi)^*\in \co*$, where * denotes monotone rearrangement.

\bD{D:4.6}
Given a collection $\mathscr P$ of substochastic matrices and an ideal $I$, we say that
$I$ is invariant under $\mathscr P$ if $(P\xi)^*\in \Sigma(I)$ for every $P\in \mathscr P$ and every $\xi \in \Sigma(I) $.
\eD

A consequence of Proposition \ref {P:4.1} is the following characterization of am-closed ideals in terms of invariance under various classes of substochastic matrices.

\bT{T:4.7}Let $I$ be an ideal. Then the following conditions are equivalent.
\item[(i)] $I$ is am-closed.
\item[(ii)] $I$ is invariant under substochastic matrices.
\item[(ii$'$)] $I$ is invariant under orthostochastic matrices.
\item[(iii)] $I\supset \mathscr L_1$ and $I$ is invariant under block orthostochastic matrices.
\eT
\bp
\item[(i)] $ \Rightarrow$ (ii) Let $\eta\in \Sigma (I)$ and let $P$ be substochastic. Then by Proposition \ref {P:4.1}, $\xi:=P\eta\in \Sigma(I^-)=\Sigma(I)$.
\item[(ii)] $ \Rightarrow$ (ii$'$) Obvious.
\item[(ii$'$)] $ \Rightarrow$ (i) Let $\xi \in \Sigma (I^-)$. By Proposition \ref {P:4.1}, $\xi:=P\eta$ for some orthostochastic matrix $P$ and some $\eta\in \Sigma (I)$. But then $P\eta\in \Sigma (I)$, which proves that $ \Sigma (I^-)\subset \Sigma (I)$, hence $ \Sigma (I^-)= \Sigma (I)$ and thus $I^-=I$, i.e., $I$ is am-closed.
\item[(ii$'$)] $ \Rightarrow$ (iii) Since (ii$'$) has just been proven above to be equivalent to (i),
$I$ is am-closed and hence contains the trace-class $\mathscr L_1$.
Being invariant under orthostochastic matrices, it is a fortiori invariant under the smaller class of block orthostochastic matrices.
\item[(iii)] $\Rightarrow$(i) Let $\xi \in \Sigma (I^-)$. By Proposition \ref{P:4.1}, $\xi:=P\eta$ for some block orthostochastic matrix $P$ and some $\eta\in \Sigma(I)$. But then $P\eta\in \Sigma (I)$, which proves that $ \Sigma (I^-)\subset \Sigma (I)$ and thus $I^- = I$, i.e., $I$ is am-closed.
\ep

Thus an ideal $I$ is am-closed if it is invariant under any class of substochastic matrices that contains the orthogonal matrices (and hence under all such classes), but invariance under a ``smaller" class may not be enough without additional hypotheses (e.g., the condition that $I\supset \mathscr L_1$ for the class of block-stochastic matrices). For instance the finite rank ideal $F$ is invariant under block stochastic matrices but is not am-closed.

It is natural and interesting to consider invariance under other important classes of substochastic matrices.
For the finite case, Birkhoff \cite {Bg46} proved that doubly stochastic matrices are convex combinations of permutation matrices.
For the infinite case, Kendall \cite {Kd60} proved that convex combinations of permutation matrices are dense in the class of doubly stochastic matrices for an appropriate topology.

It is thus natural to consider the following proper subclass of the doubly stochastic matrices (e.g., see A. Neumann's remark \cite[pg 448]{Na99}),
namely, the infinite convex combinations of permutation matrices:
\[
\mathscr C:=\{\sum_{j=1}^\infty t_j\Pi_j \mid \Pi_j \text { permutation matrix, } 0\le t_j\le1,~ \sum_{j=1}^\infty t_j=1\}.
\]

Infinite convexity of an ideal, i.e., invariance under $\mathscr C$, is in general not sufficient to guarantee that the ideal is am-closed.
For instance, it is clear that all Banach ideals are invariant under $\mathscr C$, but Banach ideals are not necessarily am-closed.
Indeed Varga \cite {Vj89} proved that for any principal ideal $(\xi)$ generated by $\diag \xi$ for some nonsummable and irregular sequence $\xi$
(we refer to \cite {Vj89} or \cite {vKgW04-Traces} for the definition),
\be{e: 18}
(\xi) \subset \cl(\xi) \subsetneq (\xi)^-= ( \cl(\xi))^-
\ee
where the closure $\cl$ is taken under the principal ideal norm which norm is well known to be complete (see \cite [Remark 3]{Vj89} and \cite[Remark 4.8] {vKgW04-Soft}). Therefore invariance under $\mathscr C$ is strictly weaker than invariance under doubly stochastic matrices.
Invariance under $\mathscr C$ is sufficient to guarantee am-closure, however, for ideals that are soft (see Definition \ref {D:soft}).

To prove this fact, we introduce the following class of block substochastic matrices:

\[
\mathscr B:=\{\sum_{k=1}^\infty \oplus t_kP_k \mid P_k \text { finite doubly stochastic matrix, } 0\le t_k\le1,~ \sum_{k=1}^\infty t_k=1\}.
\]

\bP{P:4.8} Let $I\ne\{0\}$ be an ideal.
\item[(i)] If $I$ is invariant under $\mathscr C$, then $I\supset \mathscr L_1$ and $I$ is invariant under $\mathscr B$.
\item[(ii)] If $I$ is soft-edged or soft-complemented and is invariant under $\mathscr B$, then $I$ is invariant under block stochastic matrices.
\eP

\bp
\item[(i)] Let $\textbf {t} :=~<t_k>~ \in (\ell^1)^*$ and without loss of generality assume $\sum_{k=1}^\infty t_k=1$.
Let $\Pi_k$ be the permutation matrix corresponding to the transposition $1 \leftrightarrow k$ and let \textbf {1}$:=\,<1,0, \ldots>$.
Since \textbf {1} $\in \Sigma(I)$, by the invariance of $I$ under $\mathscr C$ we have $\sum_{j=1}^\infty t_k\Pi_k \textbf {1} = \textbf {t} \in \Sigma(I)$.
This proves that $\Sigma(I)\supset (\ell^1)^* $ and hence $I\supset \mathscr L_1$.

Now let $P\in \mathscr B$, i.e., $P:= \sum_{k=1}^\infty \oplus t_kP_k$ where $0< t_k\le1,~ \sum_{k=1}^\infty t_k=1$, $P_k$ are finite doubly stochastic matrices and the k-th direct sum block corresponds to the indices $n_{k-1} < i,j\le n_k$. Recall that by the Birkhoff Theorem \cite {Bg46} (see also \cite {MO79}), each $P_k = \sum_{h=1}^{m_k} s_{h,k} \pi_{h,k}$ is a convex combination of (finite) permutation matrices $\pi_{h,k}$.
Let $\Pi_{h,k}:=I\oplus \pi_{h,k}\oplus I$ be the infinite permutation matrix corresponding to $\pi_{h,k}$, i.e., the permutation matrix that agrees with $\pi_{h,k}$ for the indices $n_{k-1}< i \le n_k$ and that leaves all the other indices invariant. For every $\xi\in \Sigma (I)$ and every $h$ and $k$,
\ba
(0\oplus \pi_{h,k}\oplus 0)\xi &= ~
< 0, 0, \ldots, 0,\pi_{h,k}\big(<\xi_{n_{k-1}+1}, \xi_{n_{k-1}+2}, \ldots, \xi_{n_k}>\big), 0, 0, \ldots> \\
&\le ~<\xi_1, \xi_2\ldots, \xi_{n_{k-1}}, \pi_{h,k}\big(<\xi_{n_{k-1}+1}, \xi_{n_{k-1}+2}, \ldots, \xi_{n_k}>\big), \xi_{n_k+1}, \ldots >
~= ~\Pi_{h,k}\xi.
\end{align*}
But then
\ba
P\xi&=\big(\sum_{k=1}^\infty \oplus t_k\sum_{h=1}^{m_k} s_{h,k} \pi_{h,k}\big)\xi \\
&=
\sum_{k=1}^\infty < 0, 0, \ldots,\big(\sum_{h=1}^{m_k} t_k s_{h,k} \pi_{h,k}<\xi_{n_{k-1}+1}, \xi_{n_{k-1}+2}, \ldots, \xi_{n_k}>\big), 0, 0, \ldots> \\
&\le \sum_{k=1}^\infty\sum_{h=1}^{m_k} t_k s_{h,k} \Pi_{h,k} \xi = \Big(\sum_{k=1}^\infty\sum_{h=1}^{m_k} t_k s_{h,k} \Pi_{h,k}\Big) \xi.
\end{align*}
Recalling that if $\eta, \zeta \in c_\mathrm o$ then $\eta \le \zeta$ implies $\eta^* \le \zeta^*$,
one now has $\sum_{k=1}^\infty\sum_{h=1}^{m_k} t_k s_{h,k} \Pi_{h,k}\in \mathscr C$ and hence
\[
(P\xi)^* \le \bigg(\big(\sum_{k=1}^\infty\sum_{h=1}^{m_k} t_k s_{h,k} \Pi_{h,k}\big) \xi \bigg)^*\in \Sigma (I),
\]
which proves that $I$ is invariant under $\mathscr B$.

\item[(ii)] Assume first that $I$ is soft-edged, let $\xi\in \Sigma(I)$ and let $P= \sum_{k=1}^\infty \oplus P_k$ be a block stochastic matrix with the k-th block supported on the indices $n_{k-1} < i,j\le n_k$.
Since $I$ is soft-edged, by definition, $\xi \le \alpha \eta$ for some $\eta\in \Sigma (I)$ and some $\alpha\in \co*$.
By grouping, if necessary, finite numbers of blocks into a single finite block and passing to a larger sequence $\alpha$ still in $\co*$,
we can assume without loss of generality that $\alpha$ is constant on each interval $(n_{k-1}, n_k]$ and that $\sum_{k=1}^\infty \alpha_{n_k} < \infty$.
By passing to a scalar multiple of $\eta$, we can furthermore assume that $\sum_{k=1}^\infty \alpha_{n_k} =1$. But then $P\xi \le P(\alpha \eta)= \big (\sum_{k=1}^\infty \oplus \alpha_{n_k}P_k\big)\eta.$
Since $R:=\sum_{k=1}^\infty \oplus \alpha_{n_k}P_k\in \mathscr B$, $(P\xi)^*\le(R\eta)^* \in \Sigma(I)$ which proves the claim.

Assume now that $I$ is soft-complemented, let $\xi\in \Sigma(I)$, and let $P= \sum_{k=1}^\infty \oplus P_k$ be a block stochastic matrix with the k-th direct sum block corresponds to the indices $n_{k-1} < i,j\le n_k$.
To prove that $(P\xi)^*\in \Sigma (I)$, it is enough to prove that $\alpha(P\xi)^*\in \Sigma (I)$ for every $\alpha\in\co*$.
The case where $\xi$ and hence $P\xi$ have finite support being trivial, assume that $\xi_n> 0$ for all $n$. Then $((P\xi)^*)_n = (P\xi)_{\Pi (n)}$ for some permutation $\Pi$.
Let $\gamma:= <\alpha_{\Pi^{-1}(n)}>$, then $\gamma \in \text{c}_o$.
Choose a sequence $\tilde{\gamma}\in\co*$ with $\gamma \le \tilde{\gamma}$ and a subsequence $n_{k_i}$ for which $r:=\sum_{i=1}^\infty \tilde{\gamma} _{n_{k_i}}< \infty$. Define $\delta_j:= \tilde{\gamma} _{n_{k_i}}$ for all $n_{k_i}< j\le n_{k_{i+1}}$ and
$R:= \sum_{i=1}^\infty \oplus \tilde{\gamma} _{n_{k_i}} \big(\sum_{h=k_i+1}^{k_{i+1}}\oplus P_h\big)$.
Then $\delta\in\co*$, $\gamma \le \delta$ and $\frac{1}{r} R\in \mathscr B$.
Since $(\gamma P\xi)_{\Pi(n)}= (\alpha(P\xi)^*)_n$ is monotone nonincreasing, it follows that $\alpha(P\xi)^*= (\gamma P\xi)^*$ and hence
$
\alpha(P\xi)^* \le ( \delta P\xi )^*= (R\xi)^* = (\frac{1}{r} R(r\xi))^*\in \Sigma(I)
$
which concludes the proof.
\ep

Thus in general, invariance under $\mathscr B$ is a weaker form of infinite convexity.

\bT{T:4.9}
If $I$ is invariant under $\mathscr C$ and is soft-edged or soft-complemented, then $I$ is am-closed.
\eT
\bp
By combining parts (i) and (ii) of Lemma \ref {P:4.8}, we see that $I$ is invariant under block stochastic matrices and that $I \supset \mathscr L_1$. By Theorem \ref {T:4.7} it follows that $I$ is am-closed.
\ep

As Varga's ideal $ \cl(\xi) $ shows (see (\ref{e: 18})), one cannot avoid the condition that $I$ is soft-edged or soft-complemented.
Indeed $\cl(\xi)$ is neither soft-edged nor soft-complemented, but being a Banach ideal it is invariant under $\mathscr C$, but by (\ref{e: 18}) it is not am-closed.

\end{document}